\documentclass[11pt]{article}
\usepackage[top=1in,bottom=1.2in,left=1.25in,right=1.25in]{geometry}
\usepackage{graphicx}
\usepackage{placeins}
\usepackage{ragged2e}
\usepackage{tabto}
\usepackage{amssymb}
\usepackage{float}
\usepackage{multirow}
\usepackage{amsmath}
\usepackage{tikz}
\usepackage{amsthm}
\usepackage{color}
\usepackage{url}
\usepackage{paralist}

\usepackage{hyperref}
\hypersetup{
    colorlinks=true,
    linkcolor=red,
    citecolor = blue,
    pdftitle={Overleaf Example},
    pdfpagemode=FullScreen,
    }
\usepackage[noabbrev,capitalise]{cleveref}
\usepackage{afterpage}
\usepackage{multicol}
\usepackage[font=footnotesize]{caption}
\usepackage[font=footnotesize]{subcaption}
\usepackage{marginnote}
\usepackage{tcolorbox}

\usepackage{newtxtext}
\usepackage{newtxmath}
\usepackage[absolute,overlay]{textpos}
\usepackage[sf,bf,medium]{titlesec}
\usepackage{bm}
\usepackage{dsfont}
\usepackage{booktabs}
\usepackage{siunitx}
\usepackage{algorithm}
\usepackage[noend]{algpseudocode}
\usepackage[textsize=tiny]{todonotes}
\usepackage[sort,compress]{cite}

\newcommand{\cosinv}{\cos^{-1}}

\newcommand{\cossq}{\cos^{2}}
\newcommand{\sinsq}{\sin^{2}}

\DeclareMathOperator\erf{Erf}

\newcommand\bLambda{\boldsymbol{\Lambda}}

\newcommand\bx{\boldsymbol{x}}

\newcommand\by{\boldsymbol{y}}

\newcommand\bF{\boldsymbol F}
\newcommand\bbf{\boldsymbol f}

\newcommand\bJ{\boldsymbol J}

\newcommand\ba{\boldsymbol a}

\newcommand\bg{\boldsymbol g}
\newcommand\bM{\boldsymbol M}
\newcommand\bU{\boldsymbol U}
\newcommand\bD{\boldsymbol D}
\newcommand\bI{\boldsymbol I}

\newcommand\bS{\boldsymbol S}
\newcommand\bPsi{\boldsymbol {\Psi}}

\newcommand\bn{\boldsymbol n}
\newcommand\br{\boldsymbol r}

\newcommand\bP{\boldsymbol P}
\newcommand\bW{\boldsymbol W}
\newcommand\bN{\boldsymbol N}

\newcommand\bX{\boldsymbol X}

\newcommand\bV{\boldsymbol V}
\newcommand\bu{\boldsymbol u}

\definecolor{darkgreen}{rgb}{0.0, 0.5, 0.13}  
\definecolor{custompink}{rgb}{0.95, 0.2, 0.6}  
\definecolor{custombrown}{rgb}{0.55, 0.16, 0.16}

\newcommand{\gb}[1]{ {\color{darkgreen} #1 }}
\newcommand{\da}[1]{ {\color{custombrown} #1 }}

\numberwithin{equation}{section}

\begin{document}


\begin{titlepage}

  \raggedleft
  {\sffamily \bfseries STATUS: {\color{red} Revised}}
  
  \hrulefill

  \raggedright
  \begin{textblock*}{\linewidth}(1.25in,2in) 
    {\LARGE \sffamily \bfseries Numerical simulation of an extensible capsule using regularized Stokes kernels and overset finite differences  \\
      \vspace{.25\baselineskip}
      }
  \end{textblock*}

  \normalsize

  \vspace{2in}
  Dhwanit Agarwal \footnote{Corresponding author at all stages. Email: dhwanit16@gmail.com; Tel.: +1 737 781 1408;} \footnote{Present address: 5 E Reed St, APT 317, San Jose, CA, 95112, United States}\\
  \emph{\small Oden Institute, University of Texas at Austin\\
    Austin, TX, United States, 78712}\\
  \texttt{\small dhwanit16@gmail.com}

   \vspace{\baselineskip}
  George Biros\\
  \emph{\small Oden Institute, University of Texas at Austin\\
    Austin, TX, United States, 78712}\\
  \texttt{\small biros@oden.utexas.edu}

  \begin{textblock*}{\linewidth}(1.25in,7in) 
    \today
  \end{textblock*}

\end{titlepage}

\begin{abstract}
In this paper, we present a novel numerical scheme for simulating deformable and extensible capsules suspended in a Stokesian fluid.  The main feature of our scheme is a partition-of-unity (POU) based representation of the surface that enables asymptotically faster computations compared to  spherical-harmonics based representations.  We use a boundary integral equation formulation to represent and discretize hydrodynamic interactions.  The boundary integrals are weakly singular. We use the quadrature scheme based on  the regularized Stokes kernels by Beale \emph{et. al.}  2019 (given in ~\cite{beale2019}). We also use partition-of unity based  finite differences that are required for the computational of interfacial forces.  Given an $N$-point surface discretization, our numerical scheme has fourth-order accuracy and $\mathcal{O}(N)$ asymptotic complexity, which is an improvement over the $\mathcal{O}(N^2 \log N)$ complexity of a spherical harmonics based spectral scheme that uses product-rule quadratures by Veerapaneni \emph{et. al.} 2011 ~\cite{veera2011}. We use GPU acceleration and demonstrate the ability of our code to simulate the complex shapes with high resolution. We study capsules that resist shear and tension and their dynamics in shear and Poiseuille flows. We demonstrate the convergence of the scheme and compare with the state of the art. 

  

\end{abstract}

\tableofcontents

\newpage

\section{Introduction}
\begin{figure}[htb!]
\centering
\includegraphics[width=0.4\linewidth]{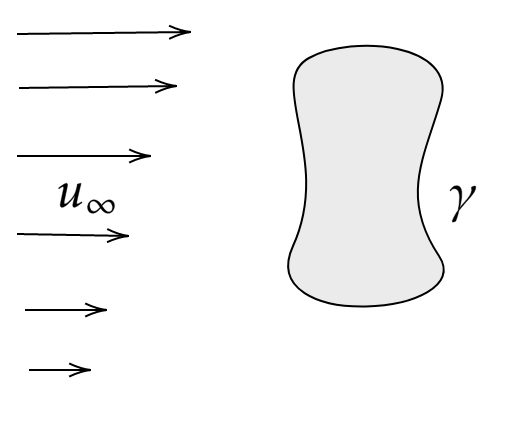}
\caption{A representation of the problem setup. The grey filled region is the interior of the capsule with membrane $\gamma$. Exterior of the capsule is filled with a Newtonian fluid and the capsule is suspended freely in it. $\bu_{\infty}$ is the imposed background fluid velocity.  }
\label{fig:schematic_domain}
\end{figure}

Capsules suspended in a Stokesian fluid  describe complex biological flows and microfluidic devices~\cite{podgorski2008, karniadakis2011,varnik2013}. In particular, we're interested in modeling red blood cell (RBC) suspensions. Several previous works~\cite{gompper2014, gekle2018, bagchi2014, fedosov2018, freund2010,rossinelli2015} have modeled  red blood cells as elastic membranes filled with a Newtonian fluid and suspended in a Newtonian fluid, also referred to as \emph{``capsules''}. In this work, we present a numerical scheme for simulating a single deformable three-dimensional capsule suspended in Stokes flow where its membrane resists shear and tension~\cite{skalak1973, poz2001} (see ~\cref{fig:schematic_domain}).

The proposed scheme allows for non-uniform discretization and fast evaluations of integral operators. We a singular quadrature based on the regularized Stokes kernel introduced in \cite{beale2019}.  We use multiple overlapping patches to parameterize the capsule surface assuming a spherical topology surface at all times. We use a fourth-order finite difference scheme using an overset grid based discretization of the patches to calculate the interfacial elastic forces. For $N$ number of discretization points of the surface, our scheme has $O(N^{2})$ work complexity, which is similar to the lower order boundary element schemes~\cite{farutin2014}. By choosing the regularization parameter in \cite{beale2019} appropriately, our overall scheme becomes fourth-order accurate.
As our singular quadrature scheme is not a product quadrature, which allows acceleration via fast multipole methods (FMMs)~\cite{greengard2008}. With FMM acceleration, our numerical scheme becomes an $O(N)$ scheme. In summary our contributions are as follows:
\begin{itemize}
    \item We describe an overlapping patch based parameterization for capsules that are  diffeomorphic to the unit sphere. 
    \item We provide a finite difference scheme based on the overset grid based discretization of these patches to calculate the surface derivatives.
    \item We evaluate  a high order singular quadrature scheme based on the regularized Stokes kernels given in ~\cite{beale2019}; combined with FMM it has an  $O(p^2)$ cost for evaluating the single layer Stokes kernel.
    \item We use GPU acceleration for the evaluation of singular quadrature to do high resolution simulations. 
\end{itemize}

\emph{Related work:} Stokesian flows with deformable capsules involve moving interfaces, fluid-structure interaction, large deformations, near and long-range stiff interactions, and often require high-order surface derivatives. Efficient methods for  simulating Stokesian capsule suspensions include immersed boundary/interface methods~\cite{bagchi2007}, lattice-Boltzmann methods~\cite{aidun2009}, and boundary integral equations methods~\cite{poz92, freund2010 ,veera2011, farutin2014}. Each of these methods trades-off aspects of  accuracy, simplicity of implementation, the ability to simulate complex physics, and numerical efficiency. Without  advocating a particular approach, we focus on  integral equation formulations and discuss the related literature in more detail. 

A nice feature of boundary integral formulations for Stokes flows is that they  require only the discretization of the capsule surface. But they require specialized singular quadrature schemes to compute layer potentials.  The methods described in~\cite{veera2011, freund2010} use spherical harmonics representations for the surface and surface fields, \emph{e.g.,}  elastic forces, velocity, and shape, to simulate the time evolution of the surface. This surface representation is spectrally accurate in the degree $p$ of the spherical harmonics that translates to  $O(p^{2})$ discretization points.
The authors ~\cite{freund2010} used the Bruno-Kunyansky quadrature~\cite{bruno2001, zorin2006} that uses a floating partition of unity to resolve weakly singular points around each global quadrature point. This scheme has $O(p^{4} \log{p})$ work complexity and can be reduced to $O(p^{3})$ using an FFT-based scheme that embeds the surface~\cite{bruno2001} in a volumetric Cartesian grid. The authors in ~\cite{veera2011} adapted the Graham-Sloan quadrature~\cite{sloan2002} to the Stokes kernel. Given $\mathcal{O}(p^2)$ global points to represent  the capsule, the scheme uses a product quadrature rule that requires $\mathcal{O}(p^2)$ spherical harmonics rotations.
\gb{For small $p$, the fastest implementation uses dense linear algebra at a $\mathcal{O}(p^6)$ cost; this becomes prohibitively expensive for large scheme. Fast rotation schemes based on Legendre transforms and FFTs or combinations of 1D non-uniform FFTS  can lower the cost to $O(p^4 \log p)$~\cite{gimbutas2013fast}.} Either way the product-quadrature has excellent accuracy, and even  for small values of $p$ for nearly-spherical shapes it outperforms the Bruno-Kunyansky quadrature. This is due to Bruno-Kunyansky's  highly localized floating-partition of unity that increases the resolution requirements. However, the product rule is too expensive for large $p$ due to the $p^4$ scaling. For example, resolving  shapes with high curvature, like the shapes shown in ~\cref{fig:red_vol60_3} and ~\cref{fig:pois_red_vol65}, requires $p>128$, which is prohibitively expensive for the product rule. The are two alternatives: change the surface representation and/or change the integration rule.

The first alternative is to avoid using global polynomials for the surface representation.  Several studies~\cite{jaeger2011, misbah2011, shaqfeh2011, seifert2001, farutin2014} have used triangulation based surface representation to provide more control over local resolution. However, these schemes  are only  second-order accurate, which  makes the computation of higher order derivatives more difficult and complicate the quadrature rules because the most accurate methods require  geometry-dependent rule precomputation and this doesn't work for time evolving geometries. 

\gb{The second approach is to change the quadrature rule. Indeed a very promising approach is the quadrature by expansion (QBX) family of methods ~\cite{wala2019fast,greengard2021fast,af2016fast}. These methods are also spectrally accurate. Direct implementations scale as $p^4$ but they can be accelerated with FMM to $p^2$. Two alternative schemes to QBX are the arbitrary-order weight-corrected trapezoidal rule \cite{wu2023unified} and the regularized Stokes kernel~\cite{beale2019}. The former is arbitrary-order accurate but a bit involved to implement. The latter that is $\sim$fourth-order accurate and relatively easy to implement.}

\emph{Limitations:}
\begin{inparaenum}[(i)]
\item
  Our scheme only works for capsules that can be smoothly mapped to the the unit sphere.  We do not consider the bending elastic forces ~\cite{poz01, farutin2014, veera2011} but our code can be extended to support them as is.
\item
  We do not consider the case of  viscosity contrast, i.e., when the fluid inside and outside the capsule have the different viscosity.  Our algorithm can be extended to include differences in the interior and exterior fluid viscosity by including a double layer Stokes potential term in the problem formulation as described in detail in \cite{poz01,rahimian2015} combined with the regularized Stokes double layer potential from \cite{beale2019}.
\item
  Our scheme has several parameters that can affect its accuracy, for example, the configuration of patches, the extent of overlap of the patches, the partition of unity functions, and a regularization parameter for the singular quadrature, whose values we set heuristically.
\item Our scheme is not curvature adaptive Instead of a regular grid, we eventually want to switch to a quadtree-based adaptive representation.
\item Our current implementation's accuracy is set by the regularized kernel. We observe the theoretical predictions for smooth shapes. For more complex shapes the asymptotic regime is not reached at the resolutions we're interested in and therefore upsampling is needed. Such an upsampling introduces a large constant in the complexity estimate.   The QBX and locally corrected trapezoidal rules can be used to deliver arbitrary accuracy and they're also compatible with FMM. The schemes would deliver optimal accuracy and complexity.
\end{inparaenum}
Alth the discussion in this paper is about a single capsule, the scheme can be extended using existing techniques to simulating an arbitrary number of capsules~\cite{veera2011} and confined boundaries~\cite{}. 

\emph{Outline of the paper:} In the next section, we state the force differential operators and the boundary integral formulation of the dynamics of an extensible capsule suspended in Stokes flow. In \cref{sec:numerics}, we describe  the surface parameterization and the numerical schemes (differentiation, integration and time stepping) we use to simulate the capsule dynamics. In ~\cref{sec:results}, we report the numerical results demonstrating the convergence and accuracy of the numerical schemes along with the results of capsule dynamics in shear and Poiseuille flow. In ~\cref{sec:FMM}, we discuss the FMM based acceleration of our scheme and show results for the speedup obtained via acceleration using a single level FMM. In ~\cref{sec:conclusion}, we provide a brief summary of the paper along with the conclusions and the future directions that can be taken to improve upon this work. 

\section{Problem formulation \label{sec:prob}}
In this section, we formally describe the mathematical formulation of the problem and introduce the notation. We discuss the differential formulation in ~\cref{sec:diff} and then specify the boundary integral formulation of the problem in ~\cref{sec:int}. A representation of the setup is shown in ~\cref{fig:schematic_domain}.  Our discussion follows the work in ~\cite{poz92, veera2011}.

\subsection{Formulation \label{sec:diff}}
Let $\gamma$ be the membrane of an extensible capsule filled with a viscous Newtonian fluid of viscosity $\mu_{i}$ and suspended in another viscous Newtonian fluid of viscosity $\mu_{e}$. In this work, we assume that the interior and exterior fluids have same viscosity, \emph{i.e.,} $\mu_{i}=\mu_{e}=\mu$. The microscopic length scale of the problem implies that the dynamics of the problem can be described by the Stokes equation. The boundary value problem is given by ~\cite{veera2011}:   

\begin{align}
        -\mu \Delta \bu(\bx) +\nabla p(\bx) &= 0\ \quad \forall \bx \in \mathbb{R}^{3} \backslash \gamma,  \label{eq:1} \\
    \nabla \cdot \bu(\bx) &= 0 \quad  \forall \bx \in \mathbb{R}^{3}\backslash \gamma \label{eq:2}, \\
    [[-p \bn +\mu  (\nabla \bu + \nabla \bu^{T})\bn]] &= \bbf  \quad \textit{ on   }  \gamma \label{eq:3}, \\
    \bu(\bx)  &\longrightarrow \bu_{\infty} \quad \textit{  for } \bx  \longrightarrow \infty \label{eq:4}, \\
    \frac{\partial \bx}{\partial t} &= \bu \quad \textit{   on } \gamma. \label{eq:5} 
\end{align}
Here $\mu$ is the viscosity of the exterior and interior fluid, $\bu$ is the velocity of the fluid, $p$ is the pressure,  [[$q$]] represents the jump in a quantity $q$ across the capsule membrane $\gamma$, $\bn$ is the unit normal to the capsule membrane, $\bbf$ is the total force exerted by the capsule membrane onto the fluid, $\bu_{\infty}$ is the imposed background velocity far away from the capsule. ~\cref{eq:1} is the Stokes equation representing the conservation of momentum, ~\cref{eq:2} is the conservation of mass, ~\cref{eq:3} is the balance of force on the interface between fluid and membrane,  ~\cref{eq:4} sets the far field velocity to be the background velocity and ~\cref{eq:5} enforces the no-slip condition on the capsule membrane. 

\emph{Interfacial force:} Now we discuss in detail the interfacial force $\bbf$ exerted by capsule membrane onto the surrounding fluid due to membrane elasticity. We only consider the elastic forces due to in-plane shear deformations of the capsule. Previous works have also included bending forces~\cite{freund2010, farutin2014, weiner78} but we leave it for future work and focus mainly on the numerical schemes in this paper. The in-plane shear force $\bbf_{s}$  is equal to the surface divergence of the symmetric part of the in-plane shear stress tensor denoted by $\Lambda$~\cite{poz2001, freund2010}. Thus, we write, 
\begin{align}
        \bbf = \bbf_{s} = \nabla_{\gamma} \cdot \Lambda \textit{.   }  \label{eq:intf} 
\end{align}
\emph{In-plane shear stress tensor $\Lambda$:} Our discussion follows the work based on  the Skalak model in ~\cite{poz92, poz2001}. The in-plane shear stress tensor $\Lambda$ is a function of the surface deformation gradient of the capsule surface relative to the stress-free reference configuration of the membrane. Let $\gamma$ be the surface in the current configuration and $\gamma_{r}$ be the reference configuration of the membrane. If $\bx_{r} \in \gamma_{r}$ is a point that maps to $\bx \in \gamma$ in the current configuration, let $\xi(\bx_{r})=\bx$ be the bijective map between the configurations. The deformation gradient $\bF$ is defined as $\bF = \frac{\partial \xi}{\partial \bx_{r}}.$
If $\bn_{r}$ is the normal to the reference configuration at $\bx_{r} \in \gamma_{r}$ and $\bn$ is the normal to the current configuration at $\bx \in \gamma$, then the relative surface deformation gradient $\bF_{S}$ is defined as $\bF_{S} = (\bI - \bn \bn^{T}) \bF (\bI - \bn_{r}\bn_{r}^{T})$, where $\bI$ is the identity tensor. It follows from the above equation that $\bF_{S}$ maps the surface tangents $\ba_{1r}$ and $\ba_{2r}$ in the reference configuration $\gamma_{r}$ to the tangents $\ba_{1}$ and $\ba_{2}$ in the current configuration $\gamma$. Also, it maps the reference normal $\bn_{r}$ to $0$. Thus, the following set of equations at every point $\bx_{r}$ uniquely determine $\bF_{S}$ at that point on the reference configuration: $\bF_{S}  \ba_{1r} = \ba_{1}, \textit{  } \bF_{S}  \ba_{2r} = \ba_{2}, \textit{  }\bF_{S}  \bn_{r} = 0.$
The left Cauchy-Green deformation tensor $\bV^{2}$ and the surface projection tensor $\bP$ at point $\bx$ in the current configuration are then defined as 
\begin{align}
    \bV^{2}(\bx) = \bF_{s}(\xi^{-1}(\bx)) (\bF_{s}(\xi^{-1}(\bx)))^{T}, \textit{  } \bP(\bx) = \bI - \bn \bn^{T}.
\end{align}
Following the work in ~\cite{skalak1973, poz2001}, if $\lambda_{1}^{2}, \lambda_{2}^{2}$ are the two non-zero eigenvalues of $\bV^{2}$, we define two scalar invariants $I_{1}= \lambda_{1}^{2}+\lambda_{2}^{2}-2 \textit{ and } I_{2}=\lambda_{1}^{2}\lambda_{2}^{2}-1$. We now define the shear stress tensor $\Lambda$ at every point $\bx$ on the current configuration as 
\begin{align}
    \Lambda(\bx) = \frac{E_{s}}{2J_{s}}(I_{1}+1)\bV^{2}(\bx) + \frac{J_{s}}{2}(E_{D}I_{2}-E_{s})\bP(\bx), \label{eq:tensor}
\end{align}
where $J_{s}=\lambda_{1}\lambda_{2}$, $E_{s}$ is the elastic shear modulus of the membrane and $E_{D}$ is the dilatation modulus. High values of the shear modulus $E_{s}$ represent higher membrane shear resistance. The dilatation modulus $E_{D}$ controls the local membrane extensibility.

\subsection{Boundary integral formulation \label{sec:int}}
Following the work in ~\cite{poz2001, poz92}, ~\cref{eq:1,eq:2,eq:3,eq:4,eq:5} can  be formulated into integral equations on the capsule membrane $\gamma$. This formulation can be written as: 
\begin{align}
   \bbf(\bx) = \nabla_{\gamma} \cdot \Lambda(\bx) \quad  \forall \bx \in \gamma  \label{eq:bi3},\\
    \bu(\bx) = \bu_{\infty}(\bx) + \mathcal{S}_{\gamma}[\bbf](\bx) \quad \forall \bx \in \gamma \label{eq:bi1}, \\
    \frac{\partial \bx}{\partial t} = \bu(\bx) \quad  \forall \bx \in \gamma \label{eq:bi2}. 
\end{align}
\da{Here $\Lambda(\bx)$ } is computed using ~\cref{eq:tensor}. $\mathcal{S}_{\gamma}[\bbf]$ represents the single layer potential of layer density $\bbf$ over the capsule membrane $\gamma$ defined as follows:
\begin{align}
    \mathcal{S}_{\gamma}[\bbf](\bx) = \int_{\gamma} \mathcal{G}(\bx,\by)\bbf(\by) d\gamma,
\end{align}
where $\mathcal{G}$ is the Stokes kernel given by $\mathcal{G}(\bx,\by) = \frac{1}{8\pi \mu}\left( \frac{I}{||\br||} + \frac{\br \otimes \br}{||\br||^{3}} \right)$ with $\br = \bx-\by$. Given the initial position of the capsule and its stress-free reference configuration, we will use  ~\cref{eq:bi3,eq:bi1,eq:bi2} to simulate the time evolution of capsule under the imposed background flow $\bu_{\infty}$. 

\section{Numerical algorithms \label{sec:numerics}}
In this section, we describe the numerical algorithms we use to simulate the dynamics of the capsule. We discretize ~\cref{eq:bi3,eq:bi1,eq:bi2} discussed in previous section. In ~\cref{sub:surfparam}, we discuss our parameterization of capsule surface (assuming it is diffeomorphic to the unit sphere) using multiple patches.  In ~\cref{sub:surfdisc}, we discuss the discretization of the surface based on our parameterization of the surface. In ~\cref{sub:nonsingint}, we discuss the numerical scheme for the surface integration of a smooth function. In ~\cref{sub:singint}, we discuss the numerical scheme for the singular integration to compute Stokes single layer potential. In ~\cref{sub:surfderiv}, we discuss the numerical scheme for calculating the surface derivatives for the computation of the elastic force.  Finally in ~\cref{sub:time_stepping}, we discuss the time stepping we use to simulate the capsule dynamics and provide a summary of the overall algorithm. 

\subsection{Surface parameterization \label{sub:surfparam}} 
We assume the capsule membrane $\gamma$ to be smooth and diffeomorphic to the unit sphere $\mathbb{S}^{2}$ embedded in $\mathbb{R}^{3}$. Thus, we have a smooth bijective  map $\phi: \mathbb{S}^{2} \longrightarrow \gamma $ with the smooth inverse $\phi^{-1}$. We define $n_{p}$ number of  overlapping patches $\{\mathcal{P}_{i}^{0}\}_{i=1}^{i=n_{p}}$, where $\mathcal{P}_{i}^{0} \subset \mathbb{S}^{2}$ form an open cover of the unit sphere $\mathbb{S}^{2}$, \emph{i.e.,} $\bigcup_{i=1}^{n_{p}}\mathcal{P}_{i}^{0} = \mathbb{S}^{2}$. Additionally, we assume each patch $\mathcal{P}_{i}^{0}$ is diffeomorphic to an open set $\mathcal{U}_{i}$ ($\subset \mathbb{R}^{2}$) via a coordinate chart $\eta_{i}^{0}$, \emph{i.e.,} $\eta_{i}^{0}:\mathcal{U}_{i} \longrightarrow \mathcal{P}_{i}^{0}$ is a diffeomorphism. The domain $\mathcal{U}_{i}$ is called the coordinate domain for the patch $\mathcal{P}^{0}_{i}$. Thus, the set $\mathcal{A}^{0} = \{(\mathcal{U}_{i}, \eta_{i}^{0})\}_{i=1}^{n_{p}}$ forms an atlas for the manifold $\mathbb{S}^{2}$~\cite{zolesio2011} . We also know such an atlas admits a smooth partition of unity subordinate to the open cover $\{\mathcal{P}_{i}^{0}\}_{i=1}^{n_{p}}$~\cite{warner83}. We choose a smooth partition of unity $\{\psi_{i}^{0}\}_{i=1}^{n_{p}}$ where
$\psi_{i}^{0}:\mathbb{S}^2 \longrightarrow \mathbb{R}$
such that $\mathsf{supp}(\psi_{i}^{0}) \subset \mathcal{P}_{i}^{0}$ for $i=1,\ldots,n_{p}$, where $\mathsf{supp}(\cdot)$ denotes the support of the function. Thus, for every $\bx_{0} \in \mathbb{S}^{2}$, $\sum_{i=1}^{n_{p}} \psi_{i}^{0}(\bx_{0}) = 1$. The precise definition of $\psi_i^0$ is given in~\cref{eq:pou0}.

\begin{figure}[t!]
\centering
\includegraphics[width=0.9\linewidth]{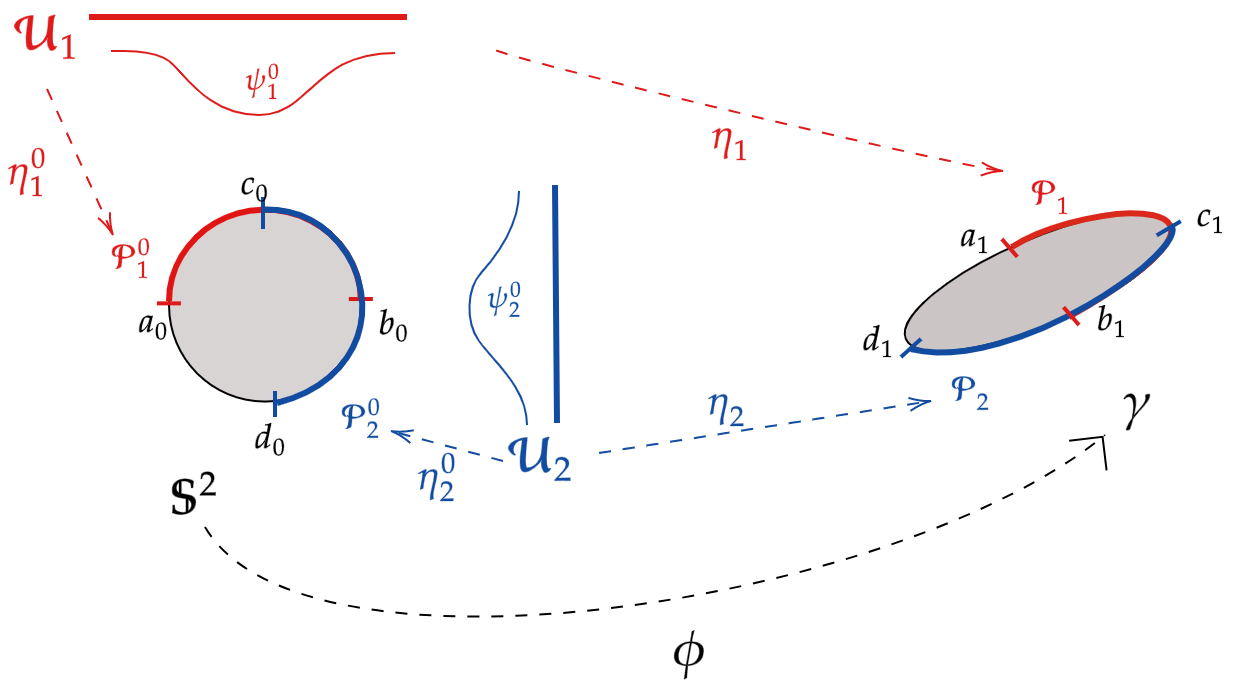}
\caption{Here we  summarize the notation for the atlas construction. The unit sphere $\mathbb{S}^{2}$  is on the left and the capsule $\gamma$ is on the right with the diffeomorphism $\phi:\mathbb{S}^{2} \longrightarrow \gamma$. We show two overlapping patches colored in red and blue. The first patch $\mathcal{P}_{1}^{0}$ is the red colored arc from point $a_{0}$ to $b_{0}$ on $\mathbb{S}^{2}$. Its corresponding patch $\mathcal{P}_{1}$ on the capsule surface $\gamma$ is shown in red as the arc from points $a_{1}$ to $b_{1}$. The second patch $\mathcal{P}_{2}^{0}$ is the blue colored arc from point $c_{0}$ to $d_{0}$ on $\mathbb{S}^{2}$. Its corresponding patch $\mathcal{P}_{2}$ on the capsule surface $\gamma$ is shown in blue as the arc from points $c_{1}$ to $d_{1}$. Their corresponding coordinate domains $\mathcal{U}_{1}$ and $\mathcal{U}_{2}$ are also shown in the red and blue color. The coordinate charts $\eta_{i}^{0}: \mathcal{U}_{1} \longrightarrow \mathcal{P}^{0}_{i}, \ i=1,2$, are shown as dashed lines in the respective colors. The diffeomorphism $\phi$ also gives the coordinate charts $\eta_{i}:\mathcal{U}_{1} \longrightarrow \mathcal{P}_{i},\  i=1,2$, for the patches on the capsule surface $\gamma$ shown as colored dashed lines from $\mathcal{U}_{i}$ to $\mathcal{P}_{i}$. The corresponding partition of unity functions $\psi^{0}_{i},i=1,2,$ with $\mathsf{supp}(\psi^{0}_{i}) \subset \mathcal{P}^{0}_{i}$ are drawn over coordinate domains $\mathcal{U}_{i}$ for visual clarity since $\mathcal{P}_{i}^{0}$ is diffeomorphic to $\mathcal{U}_{i}$.  }
\label{fig:patch_param}
\end{figure}

Define $ \mathcal{P}_{i}:= \phi(\mathcal{P}_{i}^{0})$. Now $\{\mathcal{P}_{i}\}_{i=1}^{n_{p}}$ form a set of overlapping patches that cover the  surface $\gamma$, \emph{i.e.,} $\bigcup_{i=1}^{n_{p}}\mathcal{P}_{i} = \gamma$. We then define coordinate maps $\eta_{i}: \mathcal{U}_{i} \longrightarrow \mathcal{P}_{i} $, where $\eta_{i} \equiv \phi \circ \eta_{i}^{0}$ and therefore, a corresponding atlas $\mathcal{A} := \{(\mathcal{U}_{i}, \eta_{i} )\}_{i=1}^{n_{p}}$ for the capsule surface $\gamma$. A representation of the parameterization is shown in ~\cref{fig:patch_param}. Since the patches are overlapping, a point $\bx \in \gamma$ can belong to multiple $\mathcal{P}_{i}$. For $\bx \in \gamma$, let us define the set $\mathcal{I}_{\bx} = \{1\leq i \leq n_{p} | \textit{  } \bx \in \mathcal{P}_{i} \}$. Thus, $\mathcal{I}_{\bx}$ is the collection of indices $i$ for which patch $\mathcal{P}_{i}$ contains $\bx$. We also define $\mathcal{P}_{ij} = \mathcal{P}_{i} \cap \mathcal{P}_{j}$ and $\mathcal{U}_{ij} = \eta_{i}^{-1}(\mathcal{P}_{ij})$ for $i,j=1, \ldots,n_{p}$.   We also define transition maps $\tau_{ij}: \mathcal{U}_{ij} \longrightarrow \mathcal{U}_{j}$ as $\tau_{ij} \equiv (\eta_{j}^{0})^{-1} \circ \eta_{i}^{0}|_{\mathcal{U}_{ij}}$.  Furthermore, the diffeomorphism $\phi$ allows us to create a partition of unity $\{\psi_{i}\}_{i=1}^{n_{p}}$ on $\gamma$ by defining it as 
\begin{align}
    \psi_{i}(\bx) = \psi_{i}^{0}(\phi^{-1}(\bx)),  \quad \forall \bx \in \gamma. \label{eq:pou}
\end{align}
For any function $f:\gamma \longrightarrow \mathbb{R}^{d}$, we can write 
\begin{align}
    f(\bx) = \sum_{i=1}^{n_{p}} f(\bx) \psi_{i}(\bx), \quad \forall  \bx \in \gamma.
\end{align}
If $f$ is smooth, then $f \psi_{i}$ is smooth and compactly supported in $\mathcal{P}_{i}$ for $i=1,\ldots,n_{p}$. We use this partition of unity representation to compute  derivatives and integrals on $\gamma$.

\begin{figure}[htb!]
\centering
\includegraphics[width=0.15\linewidth]{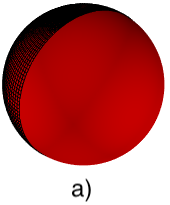}
\includegraphics[width=0.15\linewidth]{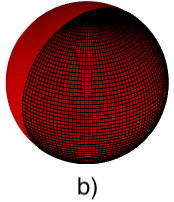}
\includegraphics[width=0.15\linewidth]{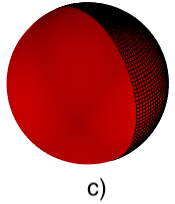}
\includegraphics[width=0.15\linewidth]{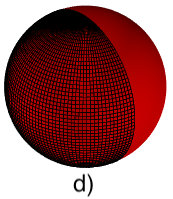}
\includegraphics[width=0.15\linewidth]{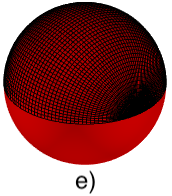}
\includegraphics[width=0.15\linewidth]{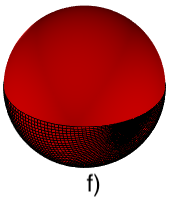}
\caption{The six hemispherical patches forming an open cover of the unit sphere $\mathcal{S}^{2}$. Each one is represented by the black grid. a)--f) $\mathcal{P}_{i}^{0}$  for   $i=1,2,\ldots,6$. }
\label{fig:patch_schematic}
\end{figure}

In this paper, we do not consider the adaptive case, and we assume that the patch parameterization remains unchanged through the calculation. In our numerical experiments, we used just six patches. Precisely, consider $\mathcal{U}_{i} = (0,\pi)\times (0,\pi), \textit{ } i=1,\ldots,6$. The coordinate maps $\{\eta_{i}^{0}\}_{i=1}^{6}$ are given below: 

\begin{align}
    \eta_{1}^{0}(u,v) &= (\sin{u} \cos{v}, \sin{u} \sin{v}, \cos{u}), \\
    \eta_{2}^{0}(u,v) &= (-\sin{u} \cos{v}, -\sin{u} \sin{v}, \cos{u}), \\
    \eta_{3}^{0}(u,v) &= (\sin{u} \sin{v}, -\sin{u} \cos{v}, \cos{u}), \\
    \eta_{4}^{0}(u,v) &= (-\sin{u} \sin{v}, \sin{u} \cos{v}, \cos{u}), \\
    \eta_{5}^{0}(u,v) &= (\sin{u} \cos{v}, -\cos{u}, \sin{u} \sin{v}), \\
    \eta_{6}^{0}(u,v) &= (\sin{u} \cos{v}, \cos{u}, -\sin{u} \sin{v}).\label{eq:param_maps}
\end{align}
These six charts form six hemispherical patches $\{\mathcal{P}_{i}^{0}\}_{i=0}^{6}$ that cover the unit sphere as shown in ~\cref{fig:patch_schematic}. The computation of the transition maps $\tau_{ij}$ is relatively simple and is given in the ~\cref{app:transition}. \da{We use the bump function, $b(r) = e^\frac{2e^{-1/|r|}}{|r|-1}$ for $|r| <1$ } and 0 otherwise, to construct the partition of unity functions on each patch. For the particular parameterization we use, for each patch of the unit sphere, we define

\begin{align} \label{eq:pou0}
    \psi_{i}^{0}(\bx_{0}) = \frac{b\left(\frac{d(\bx_{0},\eta_{i}^{0}(\pi/2,\pi/2))}{r_{0}}\right)}{\sum_{j=1}^{n_{p}} b\left(\frac{d(\bx_{0},\eta_{j}^{0}(\pi/2,\pi/2))}{r_{0}}\right)},\ \forall \bx_{0} \in \mathbb{S}^{2}
\end{align}
where $d(\bx_{0},\by_{0})$ denotes the great circle distance between  $\bx_{0}$ and $\by_{0}$ on the unit sphere, and $r_{0}>0$ determines the support of the partition of unity inside the patch. The argument of $b(\cdot)$ is  the normalized  great circle distance from a point $\bx$ to the center of the patch $\eta_{i}^{0}(\pi/2,\pi/2)$, where the normalizing factor $r_{0}$  is chosen to be $r_{0} = 5\pi/12$ (see ~\cref{app:err_pou_overlap}). Now, for a given surface $\gamma$ and diffeomorphism $\phi$, we readily have an atlas and the transition maps for the parameterization of $\gamma$. The corresponding partition of unity on $\gamma$ is available using ~\cref{eq:pou} and~\cref{eq:pou0}.


\subsection{Surface discretization \label{sub:surfdisc}} 
\begin{figure}[!t]
  \centering
\begin{subfigure}{.30\textwidth}
  \centering
  \includegraphics[width=.8\linewidth]{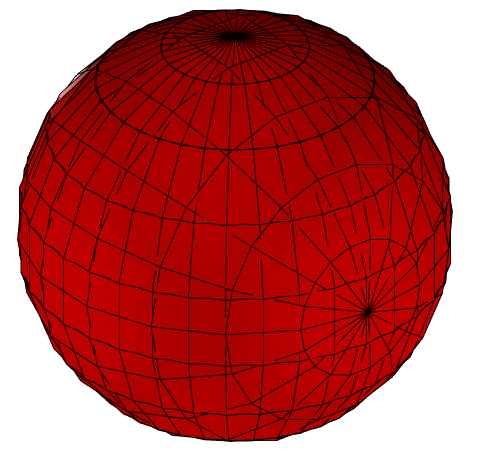}
  \caption{$m=8$}
  \label{fig:sph_mesh1}
\end{subfigure}%
\quad 
\begin{subfigure}{0.30\textwidth}
  \centering
  \includegraphics[width=.8\linewidth]{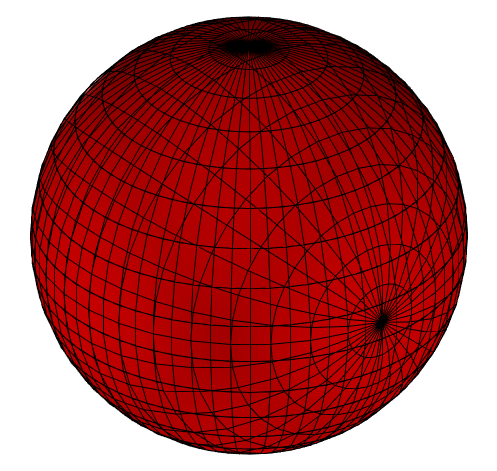}
  \caption{$m=16$}
  \label{fig:sph_mesh2}
\end{subfigure}
\begin{subfigure}{.30\textwidth}
  \centering
  \includegraphics[width=.8\linewidth]{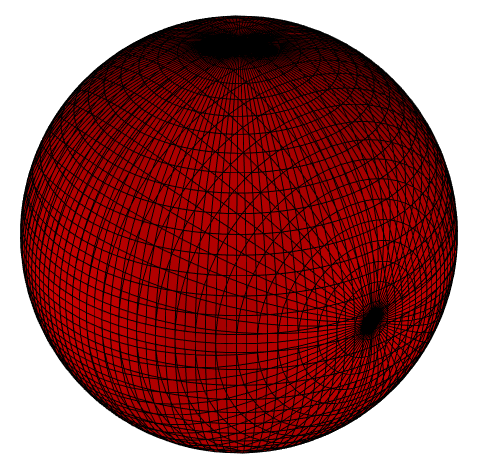}
  \caption{$m=32$  }
  \label{fig:sph_mesh3}
\end{subfigure}
\quad

\caption{Representation of discretization of a unit sphere using $m_{\mathrm{th}}$-order grids for (a) $m=8$, (b) $m=16$, (c) $m=32$. }
\label{fig:sph_mesh}
\end{figure}

\begin{figure}[!t]
  \centering
\begin{subfigure}{.30\textwidth}
  \centering
  \includegraphics[width=.8\linewidth]{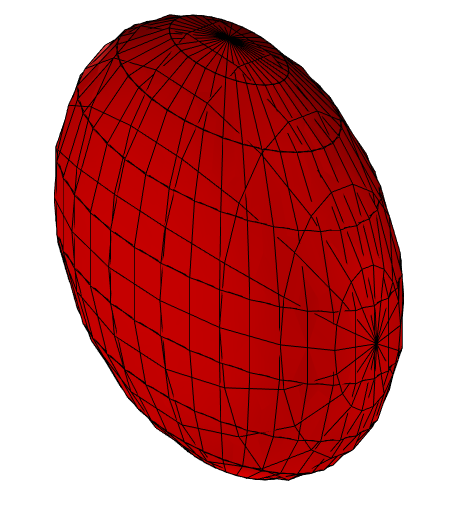}
  \caption{$m=8$}
  \label{fig:ellip_mesh1}
\end{subfigure}%
\quad 
\begin{subfigure}{0.30\textwidth}
  \centering
  \includegraphics[width=.8\linewidth]{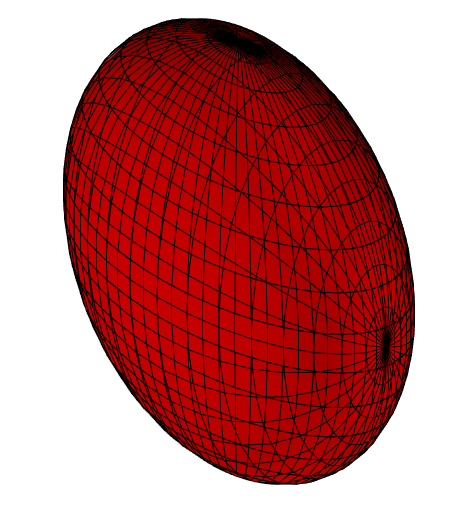}
  \caption{$m=16$}
  \label{fig:ellip_mesh2}
\end{subfigure}
\begin{subfigure}{.30\textwidth}
  \centering
  \includegraphics[width=.8\linewidth]{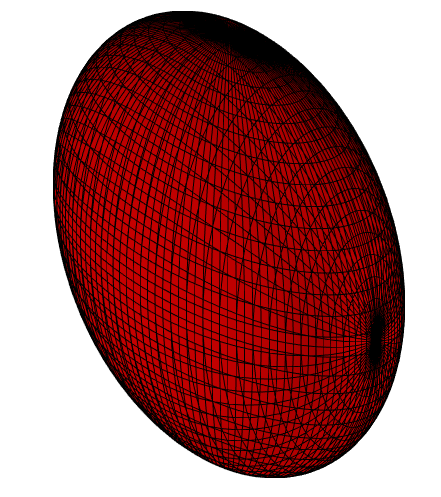}
  \caption{$m=32$  }
  \label{fig:ellip_mesh3}
\end{subfigure}
\quad

\caption{Representation of discretization of the ellipsoid $x^{2}/a^{2} + y^{2}/b^{2} + z^{2}/c^{2}=1$ with $a=0.5,b=1,c=1$, using $m_{\mathrm{th}}$-order grids for (a) $m=8$, (b) $m=16$, (c) $m=32$. }
\label{fig:ellip_mesh}
\end{figure}

We now describe in detail the discretization of the surface $\gamma$ using the parameterization  described above in ~\cref{sub:surfparam}. As mentioned above, $\mathcal{U}_{i} = (0,\pi)\times (0,\pi)$ for our parameterization.  We use the following $m_{\mathrm{th}}$-order uniform grid in $\mathcal{U}_{i}$ as follows: 
\begin{align}
    U^{m,i}_{j,k} = \left(\frac{j \pi}{m},\frac{k \pi}{m}\right),\  \forall j,k \in \{1,\ldots,m-1\}.
\end{align}
We use $h_{m}$ to denote the $m_{\mathrm{th}}$-order grid spacing in each $\mathcal{U}_{i}$ space, \emph{i.e.,} $h_{m}=\frac{\pi}{m}$. The discretization points for each patch $\mathcal{P}_{i}^{0}$ on the unit sphere are given by 
\begin{align}
    X^{0,m,i}_{j,k} = \eta_{i}^{0}(U^{m,i}_{j,k}),\  \forall j,k \in \{1,\ldots,m-1\}.\label{eq:sph_grid}
\end{align}
The discretization points for each patch $\mathcal{P}_{i}$ on $\gamma$ are given by 
\begin{align}
    X^{m,i}_{j,k} = \eta_{i}(U^{m,i}_{j,k}),\  \forall j,k \in \{1,\ldots,m-1\}.\label{eq:grid}
\end{align}
Thus, each patch contains $(m-1)^{2}$ discretization points for an $m_{\mathrm{th}}$-order grid leading to a total of $N = 6(m-1)^{2}$ points for the surface. The dynamics of the capsule are represented as the time trajectories $X^{m,i}_{j,k}(t)$. The partition of unity values at the discretization points $\psi^{m,i}_{j,k} = \psi_{i}(X^{m,i}_{j,k})$ can be precomputed and stored to be used later for computing integrals and derivatives (see \cref{eq:pou}).
The sample grids for a unit sphere and an ellipsoid are given in ~\cref{fig:sph_mesh} and ~\cref{fig:ellip_mesh}. To get the diffeomorphism $\phi$ for an ellipsoid $\gamma$ given by the level surface $\frac{x^{2}}{a^{2}} + \frac{y^{2}}{b^{2}} + \frac{z^{2}}{c^{2}} = 1$, we consider the spherical angles parameterization $\beta(u,v):[0,\pi] \times [0,2\pi) \longrightarrow \gamma$  given by
\begin{align}
    \beta(u,v) = (a\sin{u}\cos{v},b\sin{u}\sin{v},c\cos{u}), u\in [0,\pi], v\in [0,2\pi). \label{eq:sphangles_ellip}
\end{align}
We also consider the spherical angles parameterization of the unit sphere $\beta_{0}(u,v):[0,\pi] \times [0,2\pi) \longrightarrow \mathbb{S}^{2}$  given by 
\begin{align}
        \beta_{0}(u,v) = (\sin{u}\cos{v},\sin{u}\sin{v},\cos{u}), u\in [0,\pi], v\in [0,2\pi). \label{eq:sphangles_sph}
\end{align}
The diffeomorphism $\phi:\mathbb{S}^{2} \longrightarrow \gamma$ is then given as 
\[ \phi(\boldsymbol{x}_{0}) = \begin{cases} 
      \beta(\beta_{0}^{-1}(\bx_{0})) & \text{ if } \boldsymbol{x}_{0} \in \mathbb{S}^{2}\backslash \{(0,0,-1),(0,0,1)\}, \\
       (0,0,-c) & \text{ if } \boldsymbol{x}_{0} = (0,0,-1), \\
       (0,0,c) & \text{ if } \boldsymbol{x}_{0} = (0,0,1). \\
     \end{cases}
\]

\emph{Remark 1:} Given an initial capsule surface $\gamma$, we only use the diffeomorphism $\phi:\mathbb{S}^{2} \longrightarrow \gamma$ to compute the atlas $\mathcal{A}$, the partition of unity $\{\psi_{i}\}_{i=1}^{n_{p}}$ for $\gamma$ and the transition maps $\tau_{ij}$ for the coordinate domains. Once computed, we initialize the surface discretization points (see ~\cref{eq:grid}), and track their time trajectories to simulate the dynamics of the capsule. We do not need to use $\phi$ afterwards for capsule dynamics simulation. We do not change the partition of unity values after initialization. Hence, we have $\psi^{m,i}_{j,k} = \psi_{i}(X^{m,i}_{j,k}(t_{0})) = \psi_{i}(X^{m,i}_{j,k}(t))$ where $t_{0}$ is the initial time. These partition of unity values are computed at initialization and stored for subsequent use. Also, the transition maps between the coordinate domains remain unchanged and are precomputed for usage later.

\emph{Notation:} For brevity, we use $\bU^{m}$ to denote the $N \times 2$ matrix containing all the $m_{\mathrm{th}}$-order grid points, \emph{i.e.,} $\bU^{m} = [\{U^{m,i}_{j,k}\}_{1\leq j,k \leq m,  1\leq i \leq n_{p}}]^{T}$. Similarly, we define  $\bX^{m} = [\{X^{m,i}_{j,k}\}_{1\leq j,k \leq m,  1\leq i \leq n_{p}}]^{T}$ to be the $N \times 3$ matrix of corresponding discretization points on $\gamma$,  $\bPsi^{m} = [\{\psi^{m,i}_{j,k}\}_{1\leq j,k \leq m,  1\leq i \leq n_{p}}]^{T}$  to be the $N \times 1$ column vector containing the values of partition of unity functions  and $\bF^{m}$ to be the $N \times 3$ matrix containing the interfacial force $\bbf$  values at discretization points of the $m_{\mathrm{th}}$-order grid. Further, we use $\bU^{m,i}$ to denote the $(m-1)^{2} \times 2$ matrix of grid points belonging to $\mathcal{U}_{i}$. Similarly, we use $\bX^{m,i},\bPsi^{m,i}, \bF^{m,i}$ to denote the values in the patch $\mathcal{P}_{i}$. For the unit sphere, we denote the corresponding discretization points as $\bX^{0,m}$. The partition of unity column vector on $\bX^{0,m}$ is the same as $\bPsi^{m}$. 

 \emph{Remark 2:} While we do not need $\phi$ for simulating capsule dynamics after the initialization of atlas, transition maps and partition of unity, we use $||\nabla_{\mathbb{S}^{2}} \phi||_{\infty}$ as a quality metric to gauge the smoothness of the capsule surface $\gamma$ during the simulations in this paper (see ~\cref{sub:time_stepping}). To compute $||\nabla_{\mathbb{S}^{2}} \phi||_{\infty}$, we will need the discretization points $\bX^{0,m}$. Hence, we store these values as well. Note that $\phi(\bX^{0,m}) = \bX^{m}$.

\subsection{Smooth surface integrals \label{sub:nonsingint}}
Let  $f:\gamma \longrightarrow \mathbb{R}^{d}$ be  a smooth function. We can write its surface integral as a sum of the integral over all the patches using the partition of unity $\{\psi_{i}\}_{i=1}^{n_{p}}$. Then, we use the product periodic trapezoidal rule using the function values at the grid points mentioned in ~\cref{eq:grid}, which gives superalgebraic convergence~\cite{bruno2001}. We write   
\begin{align}
    \int_{\gamma} f d\gamma =  \sum_{i=1}^{n_{p}} \int_{\mathcal{P}_{i}} f\psi_{i} d\gamma =  \sum_{i=1}^{n_{p}} \int_{\mathcal{U}_{i}} f\psi_{i} W_{i} d\mathcal{U}_{i}  \approx \sum_{i=1}^{n_{p}}\left( \sum_{j,k \in \{1,\ldots,m-1\}}f(X^{m,i}_{j,k}) \psi^{m,i}_{j,k} W_{i}(X_{j,k}^{m,i})h_{m}^{2}\right),
 \label{eq:smoothint}
\end{align}
where $W_{i}$ is the surface area element for patch $\mathcal{P}_{i}$ (numerical computation of $W_{i}$ is discussed in ~\cref{sub:surfderiv}).

\emph{Work complexity:} The work complexity to compute this integral numerically  for $m_{\mathrm{th}}$-order grid is $O(n_{p}m^{2}).$ In our specific parameterization, $n_{p} = 6$. Hence, the work required is $O(6m^{2})$. 

\subsection{Singular integration \label{sub:singint}}
We use the regularized Stokes kernel described in ~\cite{beale2019} to evaluate Stokes potentials to high-order accuracy using the discretization described above in ~\cref{sub:surfdisc}. 

Consider a smooth vector-valued function $\bbf:\gamma \longrightarrow \mathbb{R}^{3}$. The single-layer Stokes potential of $\bbf$ at $\bx \in \gamma$ as described in ~\cref{sec:int} can be written as 
\begin{align}
    \mathcal{S}_{\gamma}[\bbf](\bx) = \frac{1}{8\pi \mu}\int_{\gamma}\left( \frac{\bbf(\by)}{r} + (\bbf(\by) \cdot (\bx-\by)) (\bx-\by) \frac{1}{r^{3}} \right) d\gamma,
\end{align}
where $r=||\bx-\by||$. The regularized version of Stokes single layer potential~\cite{beale2019} is given as 
\begin{align}
    \mathcal{S}_{\gamma}^{\delta}[\bbf](\bx) = \frac{1}{8\pi \mu}\int_{\gamma}\left( \frac{\bbf(\by) s_{1}(r/\delta)}{r} + (\bbf(\by) \cdot (\bx-\by)) (\bx-\by) \frac{s_{2}(r/\delta)}{r^{3}} \right) d\gamma, \label{eq:regstokes}
\end{align}
where $\delta > 0$ is the regularization parameter and $s_{1},s_{2}$ are smoothing factors such that $\frac{s_{1}(r/\delta)}{r}$ and $\frac{s_{2}(r/\delta)}{r^{3}} $ are smooth functions as $r \longrightarrow 0$. These smoothing factors ensure that $S_{\gamma}^{\delta}[\bbf]$ is an integral involving a smooth kernel and can be evaluated using the periodic trapezoidal rule as in ~\cref{sub:nonsingint}.  Following ~\cite{beale2019, cortez2014}, we choose the smoothing factors $s_{1}$ and $s_{2}$ to be 
\begin{align}
    s_{1}(r) &= \erf(r) - \frac{2}{3}r(2r^{2}-5)\frac{e^{-r^{2}}}{\sqrt{\pi}}, \\
    s_{2}(r) &= \erf(r) - \frac{2}{3}r(4r^{4}-14r^{2}+3)\frac{e^{-r^{2}}}{\sqrt{\pi}}, 
\end{align}
where $\erf(\cdot)$ is the error function. Using Taylor expansion, we can show that as $r \longrightarrow 0$, $\frac{s_{1}(r/\delta)}{r} \longrightarrow \frac{16}{3 \delta \sqrt{\pi}}$ and $\frac{s_{2}(r/\delta)}{r^{3}} \longrightarrow \frac{32}{3 \delta^{3} \sqrt{\pi}}$. This choice of smoothing factors ensures that the regularized Stokes potential is $O(\delta^{5})$ accurate~\cite{beale2019}. This regularized integral is discretized by the trapezoidal rule for smooth functions. However, the question is what should be the relation between $m$ and $\delta$? 

\emph{Regularization parameter $\delta$:} The choice of regularization parameter $\delta$ is crucial. As discussed in ~\cite{beale2016, beale2019}, the chosen $\delta$ should be large enough that regularization error dominates the discretization error in computing the integral in~\cref{eq:regstokes}. In our simulations, when computing the Stokes potential for a target point $X^{m,i}_{j,k}$ in patch $\mathcal{P}_{i}$, we choose a $\mathcal{P}_{i}$-dependent regularization parameter $\delta_{i}$ at every time step as the capsule evolves. It is defined by 
\begin{align}
    \delta^{m,i}_{j,k} = C \delta^{*} \mbox{~where~} \delta^{*} =  \left(\max_{l,l' \in \{1,\ldots,m-1\}} \{ d_{e}(X^{m,i}_{l,l'}, \mathcal{X}(X^{m,i}_{l,l'})\}\right). \label{eq:regparam} 
\end{align}
Here, we choose the constant $C=1$ (see ~\cref{app:sensitivity} for more details on the choice $C$), $\mathcal{X}(X^{m,i}_{l,l'}) =\{X^{m,i}_{l-a,l'-b} : a,b \in \{-1,0,1\} \} $ refers to the set of points which are adjacent neighbors of $X^{m,i}_{l,l'}$ and $d_{e}(a,B)$ denotes the maximum Euclidean distance between $a$ and the points in set $B$. For brevity, we define $\bm\delta^{m} = [\{\delta^{m,i}_{j,k}\}_{1\leq j,k \leq m-1,  1\leq i \leq n_{p}}]^{T}$. We experimentally observed that our choice for $\delta$ improves the accuracy the singular quadrature as in our simulations as capsule changes shapes and patches evolve over time, as opposed to using a fixed global regularization parameter $\delta = C' h_{m}$ where $C'$ is a positive constant (like in ~\cite{beale2019}). We tabulate the relative errors for different values of $C'$ in using fixed global regularization parameter $\delta = C' h_{m}$ on a sample ellipsoidal shape in ~\cref{tab:sing_int_delta} to illustrate this. 

\emph{Upsampling:} In our experiments, the observed prefactor in our weakly singular quadrature scheme is quite significant. For example at small $m$, say $m=16$, we can resolve the surface well but we need more quadrature points to resolve the weakly singular integrals. For this reason we use upsampling, especially in the small-$m$ regime.   Using numerical experiments, we determined that a four times upsampling works well. Thus, we use an upsampled grid with $N_{\mathrm{up}}=6(4m-1)^{2}$ grid points for evaluating the single layer potential.  We use cubic spline interpolation~\cite{boor78} on each patch to upsample the surface discretization points $\bX^{m}$ and surface interfacial force $\bF^{m}$. We define this upsampling operator as $\bI_{u}^{m}$. Thus, the upsampled surface points and surface force vector can be written as 
\begin{align}
    \bX^{m}_{u} = \bI_{u}^{m}  \bX^{m}, \bF^{m}_{u} = \bI_{u}^{m} \bF^{m}.
\end{align}
The partition of unity values on the upsampled grid, denoted by $\bPsi^{m}_{u}$, can also be precomputed and stored for use.  We use the upsampled values to compute the Stokes single layer potential (~\cref{eq:regstokes}) and then downsample the Stokes potential to get the values on $m_{\mathrm{th}}$-order grid. Note that both $\bI_{u}^{m}$ and $\bI_{d}^{m}$ are cubic spline based weighted interpolation operators. Their matrices can be precomputed and stored for use throughout the simulation. We also compute the $\bm \delta^{m}_{u}$ for the upsampled grid using ~\cref{eq:regparam}.  

\emph{Work complexity:} The work complexity for the evaluation of Stokes single layer potential on upsampled grid for a target point is $O(96m^{2})$. For $N_{up}$ target points in upsampled grid, the total work complexity for computing single layer potential is $ O(9216m^{4})$. The complexity for upsampling and downsampling is $ O(m^{2})$ since $\bI_{u}^{m}$ and $\bI_{d}^{m}$ are sparse linear operators. Hence, the total work complexity is $O(m^{4})$.   

\emph{GPU acceleration:} The evaluation of single layer potential on upsampled grid is the most expensive part of our numerical scheme. In order to accelerate this computation, we use a CUDA implementation of computation of Stokes single layer potential computation based on the all-pairs $O(N^{2})$ calculation using CUDA on GPU (for details refer to ~\cite{nylons2007}).

\subsection{Surface derivatives \label{sub:surfderiv}} 
In this section we discuss in detail our numerical scheme for calculating  surface derivatives. The derivatives are needed for  computing the interfacial force $\bbf$ and the surface area element $W$ (for ~\cref{eq:smoothint}) at each discretization point. The first step in the computation of the interfacial force $\bbf$ is the computation of the shear stress tensor $\Lambda$. The computation of $\Lambda$ requires the surface tangents and normals in the current configuration (see ~\cref{eq:tensor}). After computing the  stress tensor $\Lambda$, the second step is to compute the interfacial force by calculating the surface divergence of $\Lambda$ as $\bbf = \nabla_{\gamma} \cdot \Lambda$. For the computation of surface area element $W$, we need  the coefficients $E, F, G$ of the first fundamental form of the surface. The formulas for the surface divergence, denoted by $\nabla_{\gamma} \cdot $, and for the surface area element $W$ along with $E,F, G$ are summarized in the ~\cref{app:formulas}. 

Let $\eta(u,v): \mathcal{U} \longrightarrow \mathcal{P} \subset \gamma$ be a  surface patch of $\gamma$ where $\eta \in \{\eta_{1},\ldots,\eta_{n_{p}}\}$. Let $\eta(u_{0},v_{0}) = \bx \in \mathcal{P}$. Below, we summarize the steps required to evaluate the required shape derivatives $\bbf(\eta(u_{0},v_{0}))$ and $W(\eta(u_{0},v_{0}))$.

\begin{enumerate}
    \item We need to compute the tangents, denoted by $\bx_{u}$ and $\bx_{v}$, and the unit normal $\bn(\bx)$ for the current and the reference configuration. . These quantities are given by $\bx_{u} = \frac{\partial \eta}{\partial u}|_{(u_{0},v_{0})}, \bx_{v} = \frac{\partial \eta}{\partial v}|_{(u_{0},v_{0})}$ and $\bn(\bx) = \frac{\bx_{u} \times \bx_{v}}{||\bx_{u} \times \bx_{v}||}$. Given these derivatives we can then compute $\Lambda(\eta(u_{0},v_{0}))$. 
    
    \item The next step is to compute the surface divergence of $\Lambda(\eta(u,v))$ at $(u_{0},v_{0})$. This involves computation of $E,F,G$ and partial $u,v$ derivatives of the components of $\Lambda(\eta(u,v))$ (see ~\cref{app:formulas}). The computation of $E,F,G$ at point $\bx$ is straightforward using $\bx_{u}$ and $\bx_{v}$. The surface area element is then computed using $W(\bx(u_{0},v_{0}))= \sqrt{EG-F^{2}}$.  
\end{enumerate}

We note here from the above summary that the computation of both $\bbf$ and $W$ boils down to  multiple computations of the $u$ and $v$ partial derivatives of functions like the coordinate chart $\eta(u,v)$ and the components of $\Lambda(\eta(u,v))$. Now given an arbitrary function $g:\gamma \longrightarrow \mathbb{R}$, we consider the parametric function $\tilde{g}^{i} \equiv g(\eta_{i}(u,v)):\mathcal{U}_{i} \longrightarrow \mathbb{R}$ for $i=1,\ldots,n_{p}$. We denote its partial $u,v$ derivatives as $\tilde{g}^{i}_{u}$ and $\tilde{g}^{i}_{v}$. We describe below the numerical scheme to compute these partial derivatives $\tilde{g}^{i}_{u}$ and $\tilde{g}^{i}_{v}$. This scheme is then sufficient to compute $\bbf$ and $W$ at the discretization grid $\bU^{m}$ using the formulas mentioned in ~\cref{app:formulas}. Before we describe the scheme in detail, we summarize the three main steps of the scheme below:

\begin{figure}[t!]
\centering
\includegraphics[width=0.9\linewidth]{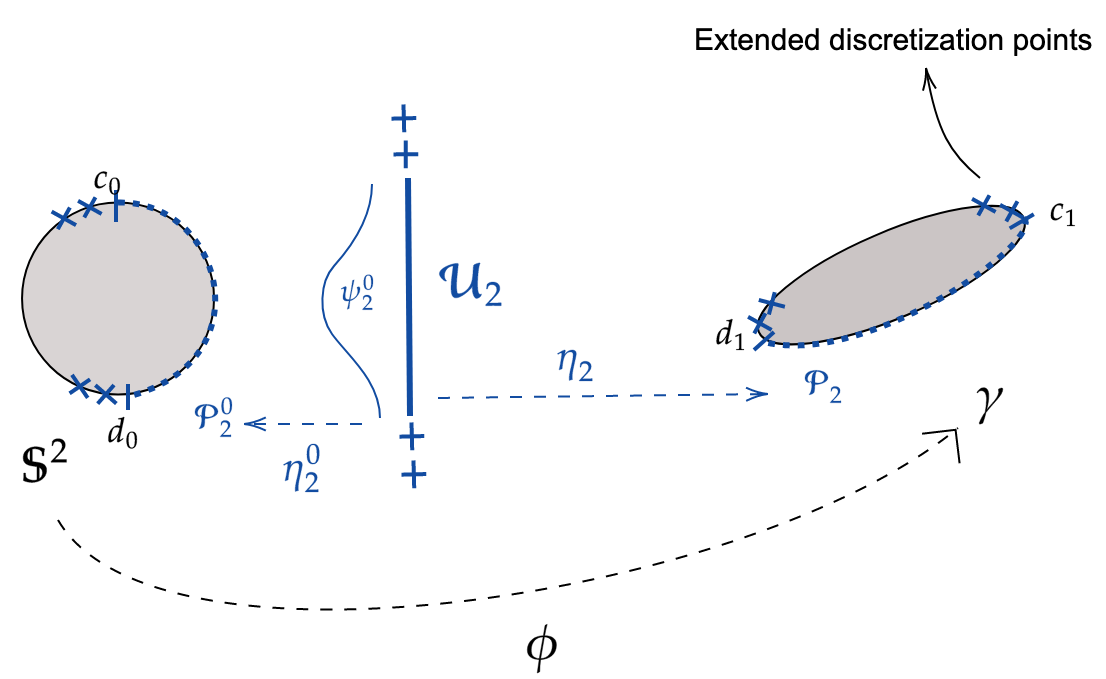}
\caption{A 2D representation of a patch extension to apply the central finite difference stencil. The unit sphere $\mathbb{S}^{2}$ is on the left and the capsule $\gamma$ is on the right with the diffeomorphism $\phi:\mathbb{S}^{2} \longrightarrow \gamma$. The discrete points on the patch $\mathcal{P}_{2}^{0}$ are represented by the blue colored dotted arc from point $c_{0}$ to $d_{0}$ on $\mathbb{S}^{2}$. Its corresponding patch $\mathcal{P}_{2}$ on the capsule surface $\gamma$ is shown in dotted blue line as the arc from points $c_{1}$ to $d_{1}$. Its corresponding coordinate domain $\mathcal{U}_{2}$ is also shown in the blue color. The extended grid points are shown in '$+$' symbol alongside $\mathcal{U}_{2}$. The extended discretization points are on $\mathbb{S}^{2}$ and $\gamma$ are shown with the '$\times$' symbol. The coordinate map $\eta_{2}^{0}$ we use in our parameterization has a natural extension to the extended domain with extended grid points with the same expression as in ~\cref{eq:param_maps}. }
\label{fig:patch_ext}
\end{figure}

\begin{enumerate}
    \item \emph{Patch extension:} We intend to use central finite difference stencil to compute $\tilde{g}^{i}_{u}$ and $\tilde{g}^{i}_{v}$ at the discretization points  for $i=1,\ldots,n_{p}$. To apply central difference at the grid points near the boundary of each coordinate domain $\mathcal{U}_{i}$, we need function values on ghost grid points outside the discretization grid in $\mathcal{U}_{i}$. We call this process of obtaining values on the extended grid points as \emph{patch extension}. We obtain the values of $\tilde{g}^{i}$ on these ghost grid points using the patch extension process. Let $\bx_{1} = \eta_{i}(u_{1}^{i},v_{1}^{i})$ for $i=1,\ldots,n_{p}$, be an extended discretization point (see ~\cref{fig:patch_ext}). The value of the function $\tilde{g}^{i}$ at the extended grid point  $(u_{1}^{i},v_{1}^{i})$, is given in the continuous form as the weighted average of the values across the patches as follows:
    \begin{align}
        \tilde{g}^{i}(u_{1}^{i},v_{1}^{i}) =  \sum_{ 1\leq j \leq n_{p}}  \tilde{g}^{j}(u_{1}^{j},v_{1}^{j}) \psi_{j}(\bx_{1}) . \label{eq:cont_ext}
    \end{align}
    
    Here, the partition of unity values are used as the weights for each patch. In the continuous case, $\tilde{g}^{i}(u_{1}^{i},v_{1}^{i}) = \tilde{g}^{j}(u_{1}^{j},v_{1}^{j})$ for $i,j=1,\ldots,n_{p}$. Hence, it is easy to see that the ~\cref{eq:cont_ext} is valid in the continuous setting by putting $\tilde{g}^{i}(u_{1}^{i},v_{1}^{i}) = \tilde{g}^{j}(u_{1}^{j},v_{1}^{j})$ for $i,j=1,\ldots,n_{p}$ $i,j=1,\ldots,n_{p}$ and noting that the partition of unity values add up to unity. In the discrete case, the above equation provides a unique definition of discretized function values at every extended grid point since $\tilde{g}^{i}(u_{1}^{i},v_{1}^{i})$ and $\tilde{g}^{j}(u_{1}^{j},v_{1}^{j})$ are not necessarily the same for $i,j\in \{1,\ldots,n_{p}\}$ in the discrete case due to numerical errors.
    
    \item \emph{Finite difference:} After extension of function values to ghost grid points, \da{we use the standard fourth order central finite difference operator} on extended grids to obtain $\tilde{g}^{i}_{u}$ and $\tilde{g}^{i}_{v}$ on the grid points $\bU^{m,i}$ for $i=1,\ldots,n_{p}$. 
    \item \emph{Blending:} We note that a point $\bx \in \gamma$ can belong to several different patches. We suppose $\bx = \eta_{i}(u_{0}^{i},v_{0}^{i})$ for $i=1,\ldots,n_{p}$. The interfacial force $\bbf(\bx)$ is intrinsic to the surface $\gamma$ and is independent of the local surface parameterization. Using patch parameterization dependent numerical derivatives $\tilde{g}^{i}_{u}$ and $\tilde{g}^{i}_{v}$ leads to different discrete derivative values at $(u_{0}^{i},v_{0}^{i})$ for different $i=1,\ldots,n_{p}$. To get unique derivative values consistent across all the patches, we \emph{blend} the derivatives  across all the patches. The blending process is taking the weighted average of derivative values across all the patches except that we also have to transform the derivatives from, say the domain $\mathcal{U}_{j}$ to the domain $\mathcal{U}_{i}$. This transformation is given by the Jacobian $\bJ_{\tau_{ij}}$ of the transition map $\tau_{ij}(u,v) = (\tau_{ij}^{(1)}(u,v), \tau_{ij}^{(2)}(u,v))  \in \mathcal{U}_{j}$ as follows: 
    \begin{align}
            \begin{bmatrix} \tilde{g}_{u}^{i} \\ 
            \tilde{g}_{v}^{i}
            \end{bmatrix} = \bJ_{\tau_{ij}} \begin{bmatrix}
                            \tilde{g}_{u}^{j} \\
                            \tilde{g}_{v}^{j}
                            \end{bmatrix}, \label{eq:jacobian_diff}
    \end{align}
    where the Jacobian matrix is
    \begin{align}
        \bJ_{\tau_{ij}} = \begin{bmatrix} \partial \tau_{ij}^{(1)}/ \partial u & \partial \tau^{(2)}_{ij}/ \partial u \\ 
                                \partial \tau_{ij}^{(1)} / \partial v & \partial \tau_{ij}^{(2)} / \partial v 
                            \end{bmatrix} .
    \end{align}
    Here, $\tau_{ij}^{(1)}, \tau_{ij}^{(2)}$ are defined as the first and second component respectively of the transition map, \emph{i.e.,} $\tau_{ij}(u,v) = (\tau_{ij}^{(1)}(u,v), \tau_{ij}^{(2)}(u,v))  \in \mathcal{U}_{j}$. Both the RHS and the LHS in ~\cref{eq:jacobian_diff} are evaluated at $(u^{i}_{0},v^{i}_{0})$. Using this Jacobian transformation, we define the blending process as the weighted average using partition of unity values as weights for the patches:
    \begin{align}
                   \begin{bmatrix} \tilde{g}_{u}^{i} \\ 
            \tilde{g}_{v}^{i}
            \end{bmatrix} = \sum_{j=1}^{n_{p}} \psi_{j}(\bx) \bJ_{\tau_{ij}} \begin{bmatrix}
                            \tilde{g}_{u}^{j} \\
                            \tilde{g}_{v}^{j}
                            \end{bmatrix}  \label{eq:jacobian_diff2}
    \end{align}
    
for all $i=1,\ldots,n_{p}$, where $\bx = \eta_{i}(u^{i}_{0},v^{i}_{0})$. The partition of unity values $\psi_{j}$ form the respective weights for the coordinate domain $\mathcal{U}_{i}$ like in the ~\cref{eq:cont_ext}. Notice that the ~\cref{eq:jacobian_diff2} is trivially valid in the continuous case since ~\cref{eq:jacobian_diff} holds and $\psi_{j}$ sum to unity.  The above equation can be expanded as follows:
 \begin{align}       
   \tilde{g}^{i}_{u}(u^{i}_{0},v^{i}_{0}) = \tilde{g}^{i}_{u}(u^{i}_{0},v^{i}_{0})
   \psi_{i}(\bx) + \sum_{j \neq i, 1\leq j \leq n_{p}}
    \left(  \Tilde{g}_{u}^{j} \left. \frac{\partial \tau_{ij}^{(1)}}{\partial u}  \right |_{(u^{i}_{0},v^{i}_{0})} +  \Tilde{g}_{v}^{j}  \left. \frac{\partial \tau_{ij}^{(2)}}{\partial u} \right \rvert_{(u^{i}_{0},v^{i}_{0})} \right)  \psi_{j}(\bx)
   \label{eq:cont_blendu}
  \\
  \tilde{g}^{i}_{v}(u^{i}_{0},v^{i}_{0}) = \tilde{g}^{i}_{v}(u^{i}_{0},v^{i}_{0}) \psi_{i}(\bx) + \sum_{j \neq i, 1\leq j \leq n_{p}} \left(  \Tilde{g}_{u}^{j} \left. \frac{\partial \tau_{ij}^{(1)}}{\partial v} \right \rvert_{(u^{i}_{0},v^{i}_{0})} +  \Tilde{g}_{v}^{j} \left. \frac{\partial \tau_{ij}^{(2)}}{\partial v} \right \rvert_{(u^{i}_{0},v^{i}_{0})} \right) \, \psi_{j}(\bx). \label{eq:cont_blendv}
\end{align}
  The terms inside the parenthesis in ~\cref{eq:cont_blendu} and ~\cref{eq:cont_blendv} come from the Jacobian transformation in ~\cref{eq:jacobian_diff}. We list all the transition maps for the specific parameterization we use in ~\cref{app:transition}. ~\cref{eq:cont_blendu} and ~\cref{eq:cont_blendv} describe the blending equation we use to get the weighted average of derivatives across all patches.  
 
\end{enumerate}

We now describe in detail all the three steps mentioned above and write the discretized versions of these steps.

\emph{Notation:} We use $\tilde{\bg}^{m}$ to denote the column vector of the values of $\tilde{g}^{i}$ at the grid points $\bU^{m,i}$  for all $i=1,\ldots,n_{p}$. We use $\tilde{\bg}^{m,i}$ to denote the column vector of the values at the grid points $\bU^{m,i}$. 

\emph{Patch extension:} As discussed above, for the discretization points near the boundary of $\mathcal{U}_{i},i=1,\ldots,n_{p}$, we need function values on \emph{ghost points}  lying outside $\mathcal{U}_{i}$ in order to use the central FDM stencil. In our simulation, we use a fifth order 7-point symmetric 1D-stencil~\cite{bengt88} for each of the partial derivatives. Thus, in order to calculate the derivative values on the grid $\bU^{m}$, we need the values of $\tilde{g}$ on the extended grid (including ghost nodes) given by 
\begin{align}
    U^{m,i,ext}_{j,k} = \left(\frac{j \pi}{m},\frac{k \pi}{m} \right),\  \forall j,k \in \{-2,-1,\ldots,m+1,m+2\},\label{eq:ext_grid}
\end{align}
 where the set of points defined by 
 $G^{m,i} := \{U^{m,i,ext}_{j,k}\}_{j,k \in \{-2,-1,\ldots,m+1,m+2\}}\backslash \{U^{m,i,ext}_{j,k}\}_{j,k \in \{1,\ldots,m-1\}} $ denotes the ghost discretization points for $\mathcal{U}_{i}$. We define the extended coordinate chart domain as $\mathcal{U}^{m,ext}_{i} = (\frac{-3\pi}{m},\frac{(m+3)\pi}{m}) \times (\frac{-3\pi}{m+1},\frac{(m+4)\pi}{m+1}) $ containing the extended grid points as defined in ~\cref{eq:ext_grid}. The parameterization of unit sphere $\eta_{i}^{0}$ has a natural extension to $\mathcal{U}_{i}^{m,ext}$ with the same function expressions as given  in \cref{eq:param_maps}. Using the diffeomorphism $\phi$, we also have mappings $\eta_{i}^{m,ext}: \mathcal{U}_{i}^{m,ext} \longrightarrow \gamma$. We define $\mathcal{P}_{i}^{m,ext}:= \eta_{i}^{m,ext}(\mathcal{U}_{i}^{m,ext})$. We define the matrix of extended grid points for all the coordinate domains as $\bU^{m,ext} = [\{U^{m,i,ext}_{j,k}\}_{-2\leq j,k \leq m+2,  1\leq i \leq n_{p}}]^{T}$, while $\bU^{m,i,ext}$ refers to the matrix of extended grid points for the coordinate domain $\mathcal{U}_{i}$. The corresponding discretization points on the surface $\gamma$ are theoretically given by
 \begin{align}
     X^{m,i,ext}_{j,k} = \eta_{i}\left(U^{m,i,ext}_{j,k}\right),\  \forall j,k \in \{-2,\ldots,m+2\}.\label{eq:ext_grid_nodes}
 \end{align}
 Let the matrix of all surface discretization nodes on extended patch be $\bX^{m,ext} = [\{X^{m,i}_{j,k}\}_{-2\leq j,k \leq m+2,  1\leq i \leq n_{p}}]^{T}$. Our goal of patch extension is to find the values of the function $\tilde{g}^{i}$ on $\bU^{m,i,ext}$, denoted by $\tilde{\bg}^{m,i,ext}$, using the known values on $\bU^{m,i}$. The discretized form of patch extension is given as 
 
 \begin{align}
         \tilde{\bg}^{m,i,ext} &= \sum_{ 1\leq j\leq n_{p}} \left(\bI_{ij}^{m,ext}  \Tilde{\bg}^{m,j}\right)  \left. \psi_{j} \right \rvert_{\bX^{m,i,ext}}.
         \label{eq:ext} 
 \end{align}
This equation is the discretized version of the ~\cref{eq:cont_ext}. Here, we additionally use the cubic splines interpolation, denoted by the operator $\bI_{ij}^{m,ext}$, to get the values of $\tilde{g}^{j}$ on $\bU^{m,i,ext}$ using the known values on the discrete grid $\bU^{m,j}$. 
 
 
 

 \emph{Finite Differences:} We use $\bD^{m,i}_{u}$ and $\bD^{m,i}_{v}$  to denote the central finite difference operators for computing partial derivative with respect to $u$ and $v$ respectively on the $m_{\mathrm{th}}$-order grid in coordinate chart domain $\mathcal{U}_{i}$. Thus, we write the partial derivatives on $\bU^{m,i}$ as 
\begin{align}
    \Tilde{\bg}_{u}^{m,i} = \bD^{m,i}_{u}  \tilde{\bg}^{m,i, ext}, \Tilde{\bg}_{v}^{m,i} = \bD^{m,i}_{v}  \tilde{\bg}^{m,i, ext}, \label{eq:fdm}
\end{align}
where $\tilde{\bg}^{m,i,ext}$ is calculated using ~\cref{eq:ext}. 

\emph{Blending:} As discussed above, applying ~\cref{eq:fdm} can lead to different numerical values of derivatives at the same discretization point for different domains $\mathcal{U}_{i}$. Therefore, we use blending to get unique values across all the patches. The discretized versions of the blending equations ~\cref{eq:cont_blendu} and ~\cref{eq:cont_blendv} are given as follows:
\begin{align}
    \tilde{\bg}_{u}^{m,i} = \Tilde{\bg}_{u}^{m,i}\left. \psi_{i}\right \rvert_{\bX^{m,i}} + \sum_{j \neq i, 1\leq j\leq n_{p}} \left( \left(\bI_{ij}^{m}  \Tilde{\bg}_{u}^{m,j}\right) \left. \frac{\partial \tau_{ij}^{(1)}}{\partial u} \right \rvert_{\bU^{m,i}} + \left(\bI_{ij}^{m}  \Tilde{\bg}_{v}^{m,j}\right) \left. \frac{\partial \tau_{ij}^{(2)}}{\partial u} \right \rvert_{\bU^{m,i}} \right) \left. \psi_{j} \right \rvert_{\bX^{m,i}}, \label{eq:blendu} \\
    \tilde{\bg}_{v}^{m,i} = \Tilde{\bg}_{v}^{m,i}\left. \psi_{i}\right \rvert_{\bX^{m,i}} + \sum_{j \neq i, 1\leq j\leq n_{p}} \left( \left(\bI_{ij}^{m}  \Tilde{\bg}_{u}^{m,j}\right) \left. \frac{\partial \tau_{ij}^{(1)}}{\partial v} \right \rvert_{\bU^{m,i}} + \left(\bI_{ij}^{m}  \Tilde{\bg}_{v}^{m,j}\right) \left. \frac{\partial \tau_{ij}^{(2)}}{\partial v} \right \rvert_{\bU^{m,i}} \right) \left. \psi_{j} \right \rvert_{\bX^{m,i}},  \label{eq:blendv} 
\end{align}
for $ i=1,\ldots,n_{p}$, where $\bI_{ij}^{m}$ is the cubic splines interpolation operator. It returns approximated function values on $\bU^{m,i}$ grid from the known values on the grid $\bU^{m,j}$. 

The three steps, namely the patch extension, the finite differences and and the blending of derivatives described above constitute our numerical scheme to calculate the interfacial forces using the formulas given in ~\cref{app:formulas}.  

\emph{Work complexity:}  The finite difference operators $\bD_{u}^{m,i}$ and $\bD_{v}^{m,i}$ and the inter-patch interpolation operators $\bI_{ij}^{m}, i,j\in \{1,\ldots,n_{p}\}$, can all be precomputed and stored for application as sparse matrices. The work complexity of patch extension is of the order of number of additional points in the extended grid, 
\emph{i.e.,} $O(n_{p}(m+5)^2 - n_{p}(m-1)^{2}) \equiv O(n_{p}m)$. The work required for applying finite difference operators is $O(n_{p}m^{2})$. Finally, the work required for blending derivatives across patches is $O(n_{p}m^{2})$. Thus, the total work complexity for calculating derivatives is $O(n_{p}m^{2})$. In our implementation $n_{p} = 6$, so the work complexity is $O(6m^{2})$.

\emph{Accuracy:} While the finite difference operators we use are fifth order accurate, our numerical scheme is expected to be fourth order accurate because we use cubic splines for patch extension and blending. This is confirmed in our numerical experiments as tabulated in ~\cref{tab:ellipsoid_centralFDM} and ~\cref{tab:4bump_centralFDM} and discussed in later sections. 

\emph{Remark 3:} We mention that the main purpose of blending is to get consistent unique derivative values by taking weighted average across different patches since a point $\bx \in \gamma$ can belong to multiple patches. We note here that the blending also improves the accuracy of the derivatives as evidenced in the results tabulated in ~\cref{tab:blending_evidence} and discussed in ~\cref{app:blending_evidence}.

\subsection{Time stepping and the overall algorithm \label{sub:time_stepping}} 
Here, we  describe the time stepping scheme we use to solve the boundary integral formulation given by ~\cref{eq:bi3,eq:bi1,eq:bi2}. The equations can be reformulated into an initial value problem as 
\begin{align}
    \frac{\partial \bx}{\partial t} &= \bu_{\infty}(\bx) + \mathcal{S}_{\gamma}[\bbf](\bx), \\
    &= \mathcal{L}(t,\bx),\  \forall \bx \in \gamma, \label{eq:bi_combined}
\end{align}
where the initial position and initial deformation gradient tensor is known. We use an explicit  Runge-Kutta-Fehlberg 4(5)~\cite{fehlberg68, fehlberg69}  adaptive time stepping scheme to simulate the time evolution of capsules. This scheme is fourth-order accurate in time. We use fixed relative tolerance $\epsilon_{tol} = 10^{-6}$ for adaptive time stepping in our simulations. Below we give the overall algorithm for the evaluation of the right hand side of ~\cref{eq:bi_combined}:

\begin{enumerate}
    \item Given the initial capsule surface, we first initialize the atlas, the transition maps and its derivatives (for \cref{eq:cont_blendu} and \cref{eq:cont_blendv}). We precompute $\bPsi^{m}, \bPsi_{u}^{m}, \bI_{u}^{m}, \bI_{d}^{m}, \bI_{i,j}^{m},  \bI_{ij}^{m,ext}, \\ \bD_{u}^{m,i}$, and $ \bD_{v}^{m,i}$ for the $m_{\mathrm{th}}$-order grid.  We also precompute the reference tangents and normals at the discretization points using the derivative scheme in ~\cref{sub:surfderiv}. Let us denote these reference tangents  by $\bX_{\partial u}^{r,m}, \bX_{\partial v}^{r,m}$ and the normals by $\bN^{r,m}$ respectively.     
    \item At a time instant $t$, we have the discretization points $\bX^{m}$ on the capsule surface . We compute the  tangents and the normals at all discretization points in the current configuration using the scheme in ~\cref{sub:surfderiv}. Let us denote these by $\bX_{\partial u}^{m}, \bX_{\partial v}^{m}$ and $\bN^{m}$. 
    \item Given the reference and the current tangents along with the normals, we compute the coefficients of the first fundamental form, the surface area element $\bW^{m}$ and the shear stress tensor $\bLambda^{m}$. 
    \item We compute the  interfacial force $\bF^{m} = \nabla_{\gamma} \cdot \bLambda^{m}$ using the derivative scheme in ~\cref{sub:surfderiv}. 
    
    \item We compute the upsampled quantities  $\bX^{m}_{u} = \bI^{m}_{u}\bX^{m}, \bF^{m}_{u}  = \bI_{u}^{m}\bF^{m}, \bW^{m}_{u} = \bI_{u}^{m}\bW^{m}$ and the regularization parameters $\bm\delta^{m}_{u}$. 
    
    \item We compute the Stokes single layer potential on the upsampled grid, denoted by $\bS^{m}_{u}$, using $\bX^{m}_{u}, \bF^{m}_{u}, \bPsi_{u}^{m}$ and  $\bW^{m}_{u}$. 
    \item We downsample the Stokes potential $\bS^{m} = \bI_{d}^{m} \bS^{m}_{u}$ and add $\bu_{\infty}(\bX^{m})$ to get the discretized version of the RHS of ~\cref{eq:bi_combined} at time instant $t$ and position $\bX^{m}$. 
\end{enumerate}

\emph{Accuracy and work complexity of the overall algorithm:} Our quadrature scheme, differentiation scheme and the time stepping scheme are all fourth order convergent. Hence, our overall scheme is fourth order convergent. We observe this convergence numerically in the convergence results in ~\cref{tab:convg_full}. The work complexity of our quadrature scheme and the differentiation scheme is $O(m^{4})$ and $O(m^{2})$ respectively for $m_{\mathrm{th}}$-order grid. Thus, the overall work complexity of our algorithm for a single time step is $O(m^{4})$. ,

\emph{Computation of $||\nabla_{\mathbb{S}^{2}} \phi||_{\infty}$}: In the results section, we monitor the norm of the surface gradient of $\phi:\mathbb{S}^{2} \longrightarrow \gamma$  with respect to the unit sphere $\mathbb{S}^{2}$ which measures the smoothness of the capsule surface during the simulation.  It serves as a proxy for measuring the stability of capsule dynamics  during the simulation. This surface gradient $\nabla_{\mathbb{S}^{2}} \phi$ is calculated using the surface gradient formula given in the ~\cref{app:formulas} by using the first fundamental form coefficients $E,F,G$ for the unit sphere and the surface area element $W$ of the unit sphere $\mathbb{S}^{2}$. We use the derivative scheme in ~\cref{sub:surfderiv} to calculate the first fundamental form using the stored discrete points $\bX^{0,m}$ on the sphere  and noting that $\phi(\bX^{0,m}) = \bX^{m}$. 




\section{Results \label{sec:results}}
Now we discuss the various results to verify the accuracy and convergence of our numerical schemes. We also present the results of simulation of extensible capsule under shear flow and the Poiseuille flow. We validate our simulations with the existing results in the literature. We also use the spherical harmonics based numerical scheme used in ~\cite{veera2011} for quantitative comparisons with our simulations. Below we give a quick summary of the spherical harmonics based scheme used in ~\cite{veera2011}.   

\emph{Spherical harmonics: }  Spherical harmonics provide an orthonormal basis~\cite{orszag74} for the square-integrable functions defined on the sphere $\mathbb{S}^{2}$. Hence, they provide a spectral representation of the surfaces and have been used for simulating Stokesian particulate flows in ~\cite{veera2011}. Below we summarize the number of discretization points used and the work complexity for the differentiation and singular integration in this scheme.
\begin{enumerate}
    \item \emph{Discretization:} We denote by $p$ the degree of spherical harmonics used to represent the surface and surface fields. For degree $p$, the scheme uses $2p(p+1)$ number of discretization points. 
    \item \emph{Singular quadrature:} The scheme uses Graham-Sloan quadrature~\cite{sloan2002} to compute the Stokes layer potential which is $O(p^{5})$ in work complexity. 
    \item \emph{Differentiation:} For surface derivatives, the scheme uses spherical harmonics based spectral differentiation which has work complexity $O(p^{3})$. Thus, the overall work complexity of the algorithm is $O(p^{5})$. 
\end{enumerate}

\subsection{Integration results}

Here, we discuss the numerical accuracy and convergence of our integration schemes discussed in ~\cref{sub:nonsingint} and ~\cref{sub:singint}. We present the relative error in computing Stokes single layer potential with and without upsampling on two different surfaces.  First, we take an ellipsoid surface given by 
\begin{align}
    \frac{x^{2}}{a^{2}} + \frac{y^{2}}{b^{2}} + \frac{z^{2}}{c^{2}} = 1, 
\end{align}
where $a=0.6, b=1, c=1$. The relative error in a quantity $q$, denoted by $\epsilon_{q}$, is defined as $\epsilon_{q} := \frac{||q - q_{\mathrm{ref}}||_{\infty}}{||q_{\mathrm{ref}}||_{\infty}}$, where $q_{\mathrm{ref}}$ is the reference solution. We use spherical harmonic solutions with degree $p=64$ as the high accuracy reference solutions for the surface derivatives~\cite{veera2011}. We report the relative errors in computing single layer Stokes potential with and without upsampling, denoted by $\epsilon^{up}_{\mathcal{S}[\bbf]}$ and $\epsilon_{\mathcal{S}[\bbf]}$ respectively, in ~\cref{tab:ellipsoid_sing}. 

\begin{figure}[t!]
\begin{subfigure}{1\textwidth}
\includegraphics[width=1\linewidth]{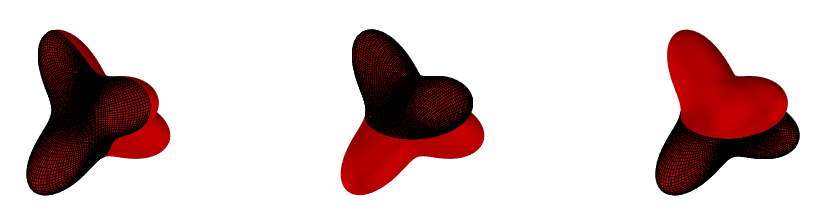}
\end{subfigure}%
\caption{The 4-bump shape with its first three patches shaded with dark mesh. }
\label{fig:flower}
\end{figure}

As a second example, we take a 4-bump shape which can be written in standard spherical parameters as 
\begin{align}
    \bX(u,v) = \rho(u,v)
\begin{pmatrix}
\sin{u} \cos{v}  \\
\sin{u}\sin{v}  \\
\cos{u} \\
\end{pmatrix}, \  \forall u \in [0,\pi] \times v \in [0,2\pi) 
\end{align}
where $\rho(u,v) = 1+e^{-3Re(Y_{lm}(u,v))}, \textit{ with } l=3,m=2$. This 4-bump shape is shown in ~\cref{fig:flower}. We report the relative errors in Stokes single layer potential on this surface in  ~\cref{tab:flo_sing}.

The upsampling improves the accuracy of the integration scheme and helps us to long time-horizon simulations that we discuss in the later sections. As evident from the tabulated results, upsampling for singular quadrature is required to get similar digits of accuracy as the derivative accuracy \footnote{The digits of accuracy for mean curvature $H$ and Gaussian curvature $K$ are lower than the upsampled quadrature but since they are not required in shear force calculation, we don't need upsampling for derivatives. However, since the bending forces~\cite{veera2011} require curvature calculations, upsampling for derivatives could be desirable if bending force is to be included. } (~\cref{tab:ellipsoid_centralFDM} and ~\cref{tab:4bump_centralFDM}).    The reference solution is computed using the Graham-Sloan quadrature~\cite{sloan2002, veera2011} for spherical harmonics order $p=64$. We also present the relative errors in computing area $A$ and volume $V$\footnote{Volume $V$ of surface $\gamma$ is given by $V = \int_{enc(\gamma)}dV$, where $enc(\gamma)$ refers to the volume enclosed by $\gamma$. Using divergence theorem, we can write it as a smooth surface integral as $V = \int_{\gamma} \bM(x,y,z)\cdot \bn d\gamma$, where $\bM(x,y,z) = (x,0,0)$.  } to demonstrate the convergence and accuracy of our integration scheme for smooth functions. For further verification and to study the change in relative errors with the reduced volume $\nu$, we provide the error in our upsampled quadrature scheme in ~\cref{app:err_redvol}. In ~\cref{app:err_shear_pois_shapes}, we also tabulate the relative errors demonstrating the convergence of our derivative and integration schemes for the complex shapes obtained in actual shear and Poiseuille flow simulations discussed in the later sections. This further verifies the correctness of our code. 

\begin{table}
\centering
\caption{Relative error in the computation of  Stokes single layer potential (with density $\bbf = (x^{2},y^{2},z^{2})$ on (a) an ellipsoid and (b) the 4-bump shape (as shown in ~\cref{fig:flower}). $\epsilon_{\mathcal{S}[\bbf]}$ is relative error with no upsampling and $\epsilon^{up}_{\mathcal{S}[\bbf]}$ is the relative error with 4$\times$ upsampling. $\epsilon_{A}$ and $\epsilon_{V}$ are the relative errors in computing area and volume respectively.  $N$ is the number of discretization points. Reference solution is computed using $p=64$ spherical harmonics (see ~\cite{veera2011}).}
\begin{tabular}{SSSSSS} \toprule
    {$m$}  &  {$N$} & {$\epsilon_{\mathcal{S}[\bbf]}$  } &  {$\epsilon^{up}_{\mathcal{S}[\bbf]}$} & {$\epsilon_{A}$} & {$\epsilon_{V}$}  \\ \midrule
    8   &  294 & 6E-2  &  7E-3 & 5E-3  & 3E-3 \\
    16  & 1350 & 8E-3 &  4E-4  & 5E-4  & 4E-4 \\    
    32  &  5766 & 4E-4 & 2E-5  &  1E-5  & 1E-5 \\
    64  & 23814 & 1E-5 &  1E-6 & 1E-6   & 1E-6   \\ \bottomrule

\end{tabular}
\subcaption{Relative error in Stokes single layer potential on ellipsoid of reduced volume $\nu=0.9$, given $\frac{x^{2}}{a^{2}} + \frac{y^{2}}{b^{2}} + \frac{z^{2}}{c^{2}} = 1$, where $a=0.6, b=1, c=1$. \label{tab:ellipsoid_sing} }
\vspace{1cm}
\begin{tabular}{SSSSSS} \toprule
    {$m$}  &  {$N$} & {$\epsilon_{\mathcal{S}[\bbf]}$  } &  {$\epsilon^{up}_{\mathcal{S}[\bbf]}$} & {$\epsilon_{A}$} & {$\epsilon_{V}$}  \\ \midrule
    8   & 294  &5E-1  &  4E-1 & 2E-1  & 2E-1   \\
    16  &  1350 & 2E-1  &  9E-2 & 3E-2  & 4E-2   \\    
    32  & 5766  &3E-2 &  8E-3  & 2E-3  &  3E-3 \\
    64  &  23814  &3E-3 &  6E-4  & 2E-4  & 2E-4   \\ \bottomrule

\end{tabular}
\subcaption{Relative error in Stokes single layer potential for the 4-bump shape. \label{tab:flo_sing}} 
\end{table}

\subsection{Derivative results \label{sub:res_deriv}}
Now we discuss the results showing the accuracy and convergence of our numerical differentiation scheme. We report the relative errors in derivative calculations on the ellipsoid and the 4-bump shape in  ~\cref{tab:ellipsoid_centralFDM} and ~\cref{tab:4bump_centralFDM} respectively. The derivatives converges as the the grid length $h_{m}$ decreases. Empirically, the convergence is of the fourth order which is what we predicted because of the cubic spline interpolation used in patch extension operators and blending. We also provide verification of derivatives for ellipsoids of different reduced volumes and tabulate the results in ~\cref{app:err_redvol}.

\begin{table}
\centering
\caption{Relative error in the derivative scheme for (a) an ellipsoid and (b) the 4-bump shape (shown in ~\cref{fig:flower}). We report relative error in surface normals $\bn$, mean curvature $H$, Gaussian curvature $K$ and surface divergence $div_{\gamma}$. $N$ is the number of discretization points. Spherical harmonics results for $p=64$, using the algorithms in ~\cite{veera2011}, are used as true values and the error is computed relative to those values. Surface divergence is computed for the smooth function $\bg(x,y,z)=(x^{2},y^{2},z^{2})$ on the surface.  }
\begin{tabular}{SSSSSS} \toprule
    {$\textit{m}$} & {$N$} & {$\epsilon_{\bn}$}  & { $\epsilon_{H}$} & {$\epsilon_{K}$} & { $\epsilon_{div_{\gamma}}$}\\ \midrule
    8    & 294  &  5E-3 &    4E-3  & 6E-3 &  4E-2 \\
    16   & 1350  & 2E-4  &   5E-4   &  1E-4 &  3E-3 \\    
    32   & 5766  &  1E-5  &    3E-5  & 6E-5  &  1E-4 \\
    64   & 23814  & 7E-7  &    2E-6   &  3E-6  &  1E-5  \\  \bottomrule

\end{tabular}
\subcaption{Relative error in derivatives for an ellipsoid of reduced volume $\nu=0.9$, given by $\frac{x^{2}}{a^{2}} + \frac{y^{2}}{b^{2}} + \frac{z^{2}}{c^{2}} = 1$, where $a=0.6, b=1, c=1$. \label{tab:ellipsoid_centralFDM}}
\vspace{1cm}
\begin{tabular}{SSSSSS} \toprule
    {$\textit{m}$} & {$N$} & {$\epsilon_{\bn}$} &  { $\epsilon_{H}$} & {$\epsilon_{K}$} & { $\epsilon_{div_{\gamma}}$}\\ \midrule
    8    &  294 & 2E-1 &   3E-1  & 8E-1 &  6E-1 \\
    16   &  1350 & 6E-2  &   1E-1   &  1E-1 &  3E-1 \\    
    32  & 5766 & 6E-3   &    1E-2  & 2E-2  &  6E-2 \\
    64   & 23814 & 5E-4  &   1E-3   &  2E-3  &  5E-3  \\  \bottomrule

\end{tabular}
\subcaption{Relative error in derivatives for the 4-bump shape. \label{tab:4bump_centralFDM}} 
\end{table}

\subsection{Accuracy and convergence of the full numerical scheme }
\begin{figure}[t!]
\begin{subfigure}{0.25\textwidth}
\includegraphics[width=1\linewidth]{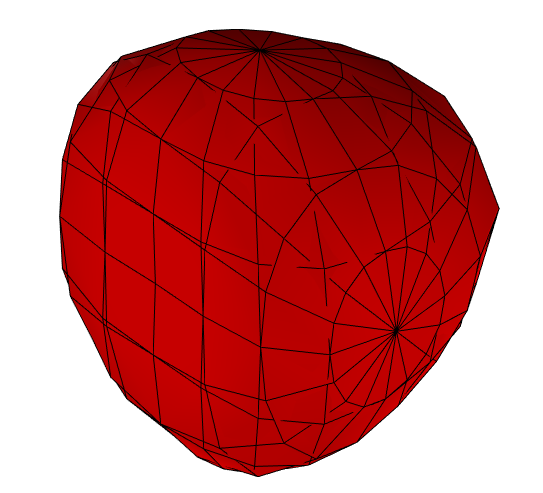}
\subcaption{}
\end{subfigure}%
\begin{subfigure}{0.25\textwidth}
\includegraphics[width=1\linewidth]{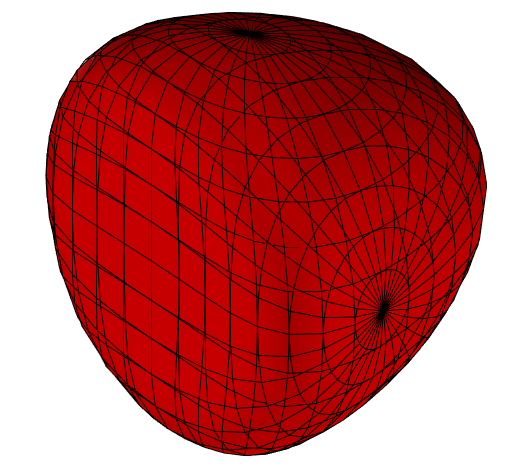}
\subcaption{}
\end{subfigure}%
\begin{subfigure}{0.25\textwidth}
\includegraphics[width=1\linewidth]{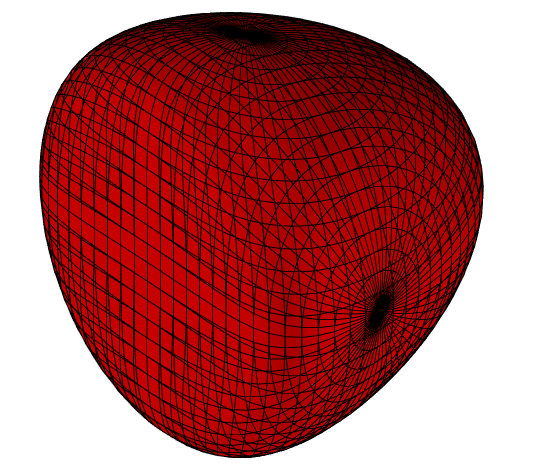}
\subcaption{}
\end{subfigure}%
\begin{subfigure}{0.25\textwidth}
\includegraphics[width=1\linewidth]{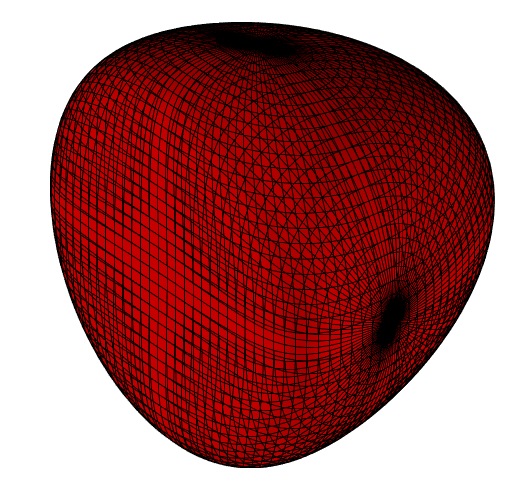}
\subcaption{}
\end{subfigure}%
\caption{The parachute shape under background Poiseuille flow (see ~\cref{eq:poiseuille_flow}) obtained after time $T=0.5$ for different levels of discretization (a) $m=8$, (b) $m=16$, (c) $m=32$ and (d) $m=48$ for $E_{s}=2, E_{D}=20$. The initial shape is same as the stress-free reference state and is given by $x^{2}/a^{2} + x^{2}/b^{2} + x^{2}/c^{2} = 1,$ where $a=0.9,b=1.0,c=1.0$. }
\label{fig:convg_poiseuille_shapes}
\end{figure}

Now, we discuss the accuracy and convergence of our full numerical solver for the extensible capsule simulation. To this end, we take an ellipsoidal capsule with an initial surface configuration given by $\frac{x^{2}}{a^{2}} + \frac{y^{2}}{b^{2}} + \frac{z^{2}}{c^{2}} = 1$, where $a=0.9, b=1, c=1$. The initial shape is taken to be the stress-free reference configuration. We do the simulations under two different imposed background flows, \emph{i.e,} shear flow and Poiseuille flow. The flows are described mathematically  as 
\begin{align}
     \bu_{\infty}(x,y,z) &= \Dot{\gamma}(y,0,0), \textit{ (shear flow) }, \label{eq:shear_flow}\\
    \bu_{\infty}(x,y,z) &= (\alpha (R_{0}^{2} - y^{2} - z^{2}),0,0) \textit{ (Poiseuille flow) }, \label{eq:poiseuille_flow}
\end{align}
where $\Dot{\gamma}$ is the shear rate of the shear flow, $\alpha$ controls the curvature of the Poiseuille flow and $R_{0}$ is the radius of the circular cross-section of the Poiseuille flow. We set $\Dot{\gamma} = 1$, $\alpha=1$ and $R_{0}=5$ for these simulations. We use the differentiation, integration and the time stepping schemes discussed above to simulate these setups for a fixed time horizon $[0,T]$ using our solver. In ~\cref{tab:convg_full}, we tabulate the relative errors in area $A$, volume $V$ and moment of inertia tensor $J$ of the final capsule shape for different values of the discretization  order $m$. We observe fourth order convergence in the relative errors. The parachute shapes at time $T=0.5$ for the Posieuille flow simulations for different levels of discretization are  given in ~\cref{fig:convg_poiseuille_shapes}.  

\begin{table}
\centering
\caption{Self-convergence results for an extensible capsule simulation under background (a) shear flow (given by $\bu_{\infty}(x,y,z) = (y,0,0)$)  and (b) Poiseuille flow (given by $\bu_{\infty}(x,y,z) = ( (25 - y^{2} - z^{2}),0,0)$).  The initial shape is an ellipsoid given by $x^{2}/a^{2} + x^{2}/b^{2} + x^{2}/c^{2} = 1,$ where $a=0.9,b=1.0,c=1.0$. We set the shear modulus to be $E_{s}=2$ and the dilatation modulus to be $E_{D}=20$. We simulate the capsule for a time horizon $[0,T]$ and report the relative errors in area $A$, volume $V$ and moment of inertia tensor $J$ of the capsule shape at time $T=0.5$. The reference solution for the relative error is the numerical solution computed using the $m=48$ grid. We use Runge-Kutta-Fehlberg time stepping scheme ( see ~\cref{sub:time_stepping}) with fixed relative tolerance of $\epsilon_{tol} = 10^{-6}$.   \label{tab:convg_full}}
\begin{tabular}{SSSS} \toprule
    {$m$}  & {$\epsilon_{A}$  } &  {$\epsilon_{V}$  } &  {$\epsilon_{J}$}  \\ \midrule
    8    &  2E-2 & 4E-2  & 3E-2   \\
    16   &  2E-3 & 1E-3  & 2E-3   \\    
    32   &  1E-4  & 1E-4  &  1E-4  \\ \bottomrule

\end{tabular}
\subcaption{ Relative errors in the capsule dynamics driven by a shear flow.  \label{tab:convg_shear} }
\vspace{1cm}
\begin{tabular}{SSSS} \toprule
    {$m$}   & {$\epsilon_{A}$  } &  {$\epsilon_{V}$  } &  {$\epsilon_{J}$}  \\ \midrule
    8    &  2E-2 & 4E-2  & 2E-2   \\
    16   &  3E-3 & 1E-3  & 2E-3   \\    
    32   &  2E-4  & 1E-4  &  1E-4  \\ \bottomrule

\end{tabular}
\subcaption{ Relative errors in  the capsule dynamics driven by a background Poiseuille flow. \label{tab:convg_poiseuille}} 
\end{table}

\subsection{Relaxation of capsule }
\begin{figure}[t!]
\begin{subfigure}{0.4\textwidth}
\includegraphics[width=1\linewidth]{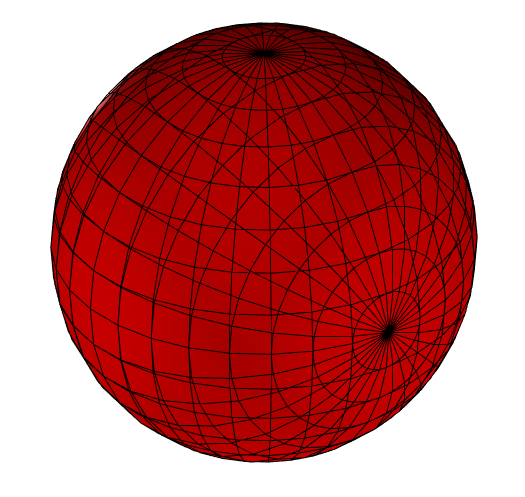}
\subcaption{$t=0$}
\end{subfigure}%
\begin{subfigure}{0.35\textwidth}
\includegraphics[width=1\linewidth]{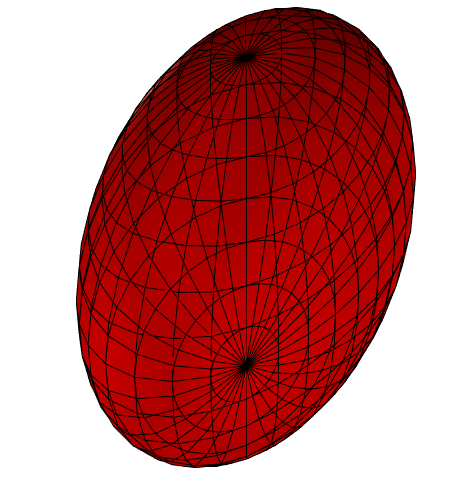}
\subcaption{$t=T_{1}$}
\end{subfigure}
\begin{subfigure}{0.4\textwidth}
\includegraphics[width=1\linewidth]{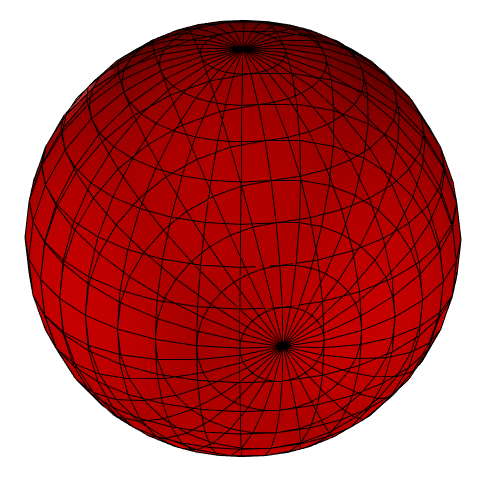}
\subcaption{$t=T$}
\end{subfigure}
\begin{subfigure}{0.6\textwidth}
\includegraphics[width=1\linewidth]{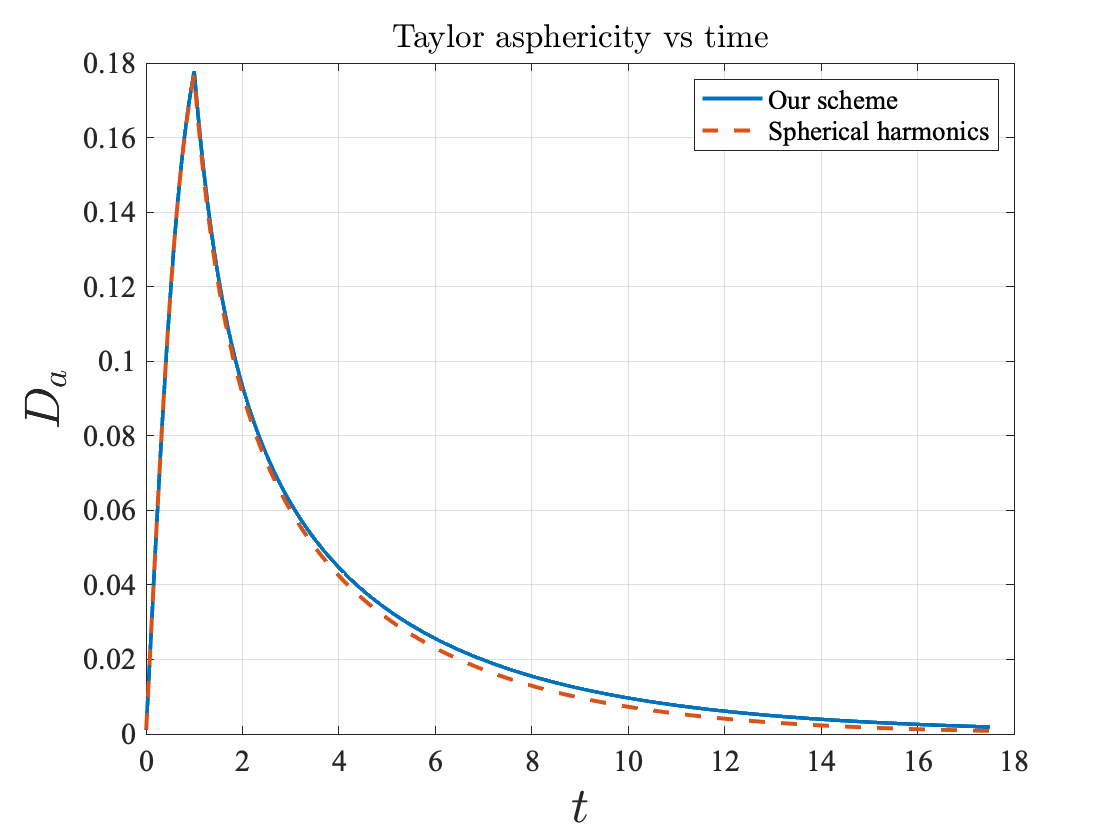}
\subcaption{ \label{fig:taylordef_relax_plot}}
\end{subfigure}
\caption{ Shapes obtained at different times during simulation of the relaxation of the capsule with initial spherical shape as the stress-free reference state. We impose background shear flow $\bu_{\infty}(x,y,z) = (y,0,0)$ for time horizon $[0,T_{1}]$. This is followed by zero background velocity from $[T_{1},T]$. The capsule relaxes back to the stress-free reference shape as expected. We take $m=16, T_{1}=1, T=17.5, E_{s}=2$ and $E_{D}=20$. (a) Shape at time $t=0$. (b) Shape at time $t=T_{1}$. (c) Shape at time $t=T$. (d) The plot of Taylor asphericity $D_{a}$ vs time $t$, where $D_{a}$ increases till $t=1$ under background shear flow and then drops to zero during the relaxation phase when the background flow is zero. Blue solid lines is the plot using our scheme and red dotted lines denote is the plot using spherical harmonics based scheme (with spherical harmonic degree $p=32$).  }
\label{fig:relaxation}
\end{figure}

Now, we simulate the relaxation of a capsule to its stress-free reference state to validate our code. We take a capsule with an initial shape of a unit sphere. The initial shape is also the stress-free reference configuration of the capsule. We simulate the dynamics of the capsule under a background linear shear flow with shear rate $\Dot{\gamma}=1$ for a fixed time horizon $[0,T_{1}]$ followed by a zero background flow velocity for $[T_{1},T]$. We report the resulting shapes at different times in ~\cref{fig:relaxation} . We observe that once the background flow diminishes the capsule returns to the stress-free state, \emph{i.e.,} a unit sphere. Using the moment of inertia tensor $J$ and volume $V$ of the capsule, we compute the instantaneous Taylor asphericity parameter~\cite{lac2004} $D_{a}$ of the capsule shape defined as 
\begin{align}
    D_{a} = \frac{L-S}{L+S}, \\
    \textit{ where } S = \sqrt{\frac{J_{xx} + J_{yy} - \sqrt{(J_{xx}-J_{yy})^{2} + 4J_{xy}^{2}}}{2V}}&, \textit{   } L = \sqrt{\frac{J_{xx} + J_{yy} + \sqrt{(J_{xx}-J_{yy})^{2} + 4J_{xy}^{2}}}{2V}}.
\end{align}
For a sphere, $D_{a}=0$. We plot the Taylor asphericity as a function of time to monitor the relaxation back to the reference unit sphere shape and plot it in ~\cref{fig:taylordef_relax_plot}. As expected, we observe that the capsule relaxes back to a unit sphere when the background flow is removed as $D_{a}$ drops back to zero. Our results agree quantitatively with the results obtained using spherical harmonics based scheme used in ~\cite{veera2011}.

\subsection{Capsule in shear flow}
\begin{figure}[t!]
\begin{subfigure}{0.55\textwidth}
\includegraphics[width=1\linewidth]{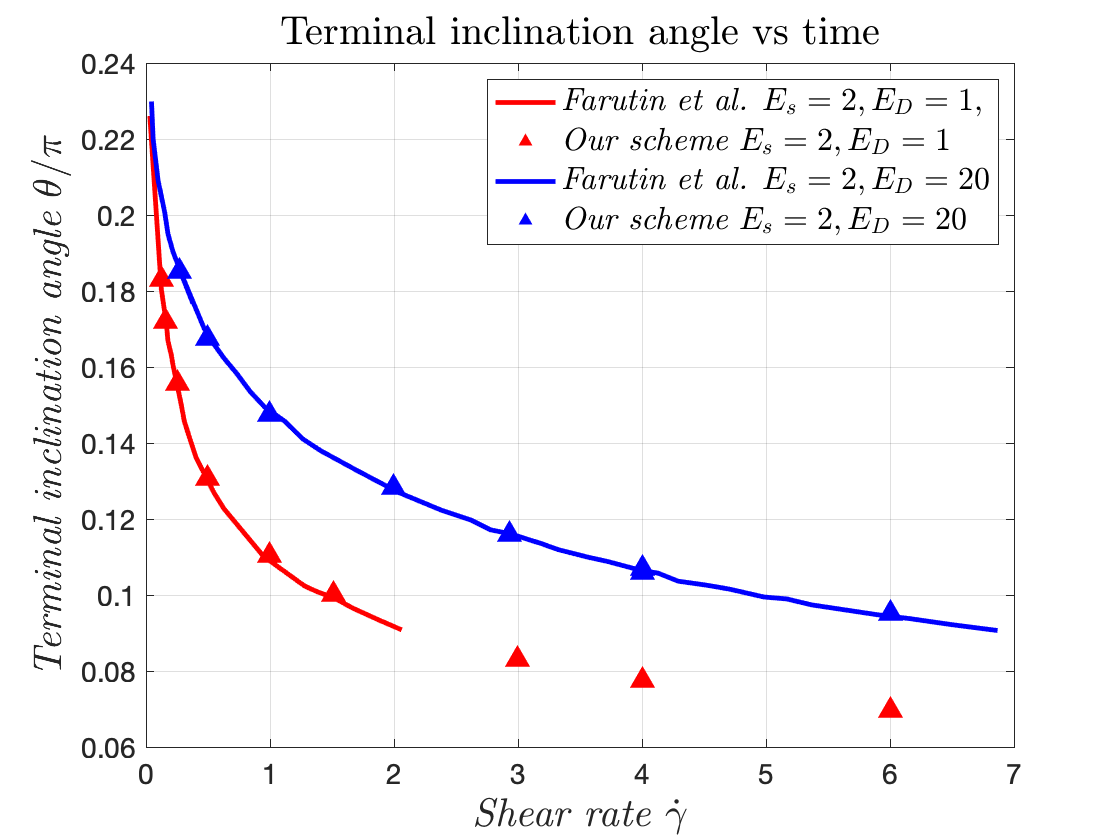}
\subcaption{\label{fig:shear_inclination}}
\end{subfigure}%
\begin{subfigure}{0.60\textwidth}
\includegraphics[width=1.13\linewidth]{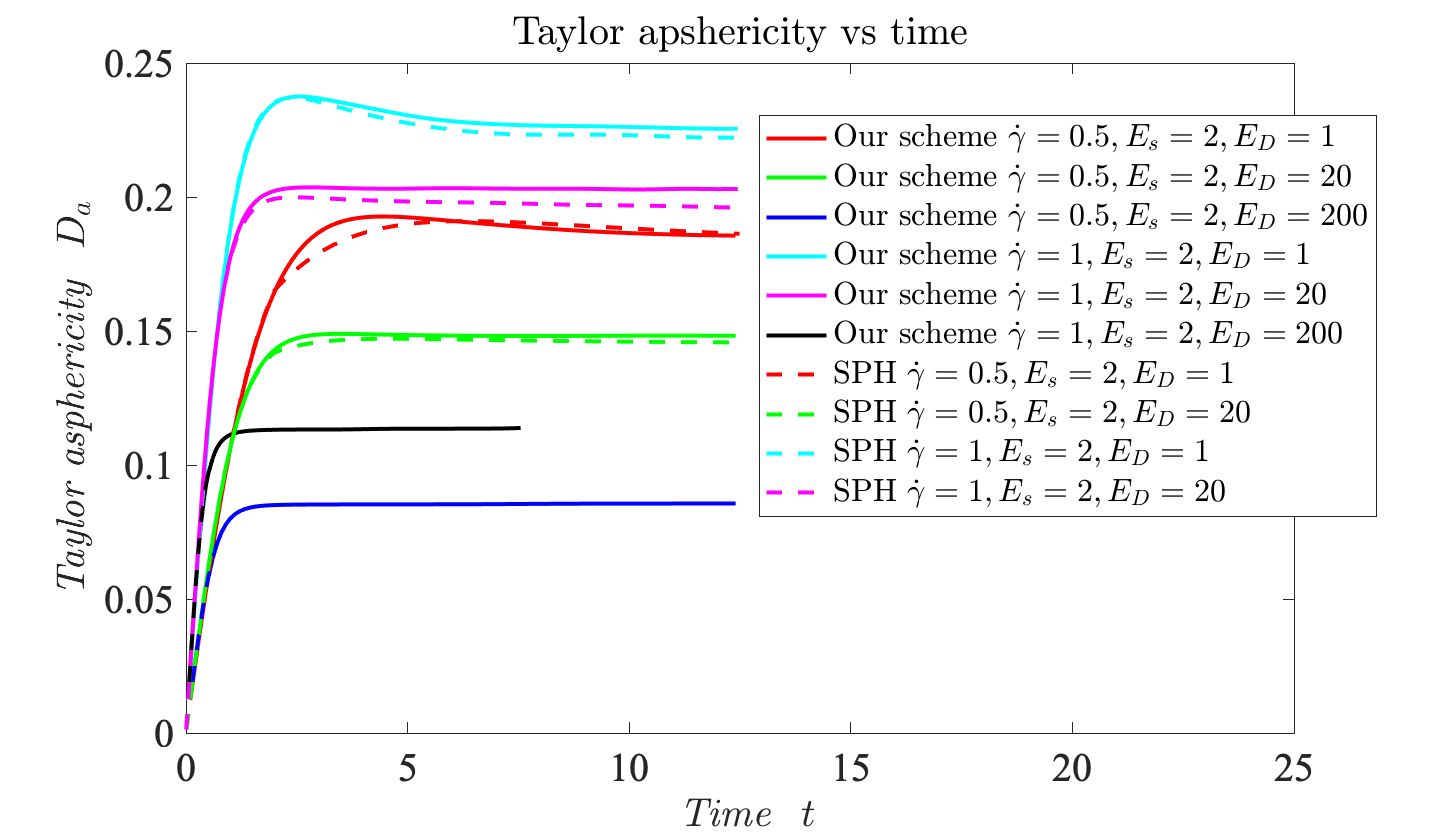}
\subcaption{\label{fig:shear_asphericity}}
\end{subfigure}
\caption{ Validation of shear flow simulation results computed using $m=16$. (a) Plot of terminal inclination angle $\theta$ for different shear rates and membrane mechanical properties ($E_{s}$ and $E_{D}$). We compare our results with results in ~\cite{farutin2014} and obtain good quantitative agreement. (b) We plot the evolution of Taylor asphericity $D_{a}$ of the the capsule vs time and compare our results (solid lines) with the spherical harmonics ($p=32$) based scheme in ~\cite{veera2011} (dashed lines). We observe good quantitative agreement. For $E_{D} =200$, the spherical harmonics code is unable to resolve the shapes and the surface fields and therefore, we do not provide plots for those.   }
\label{fig:shear}
\end{figure}

\begin{figure}[t!]
\begin{subfigure}{0.30\textwidth}
\includegraphics[width=1\linewidth]{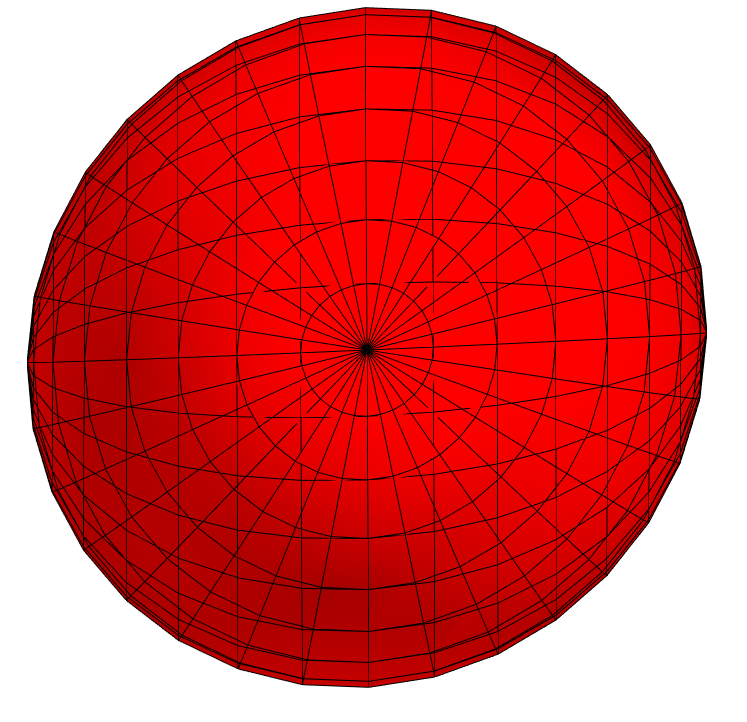}
\subcaption{$t=0$\label{fig:red_vol96_1}}
\end{subfigure}
\begin{subfigure}{0.30\textwidth}
\includegraphics[width=1\linewidth]{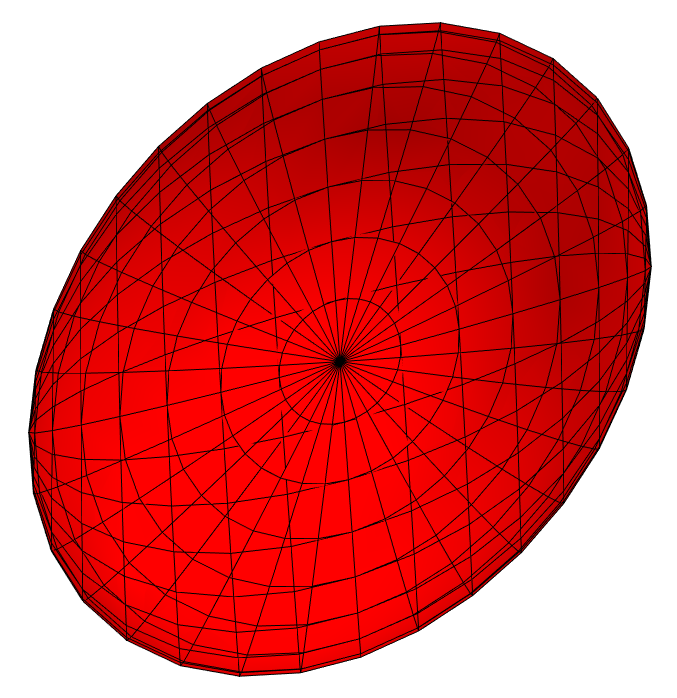}
\subcaption{$t=0.25$\label{fig:red_vol96_2}}
\end{subfigure}
\begin{subfigure}{0.35\textwidth}
\includegraphics[width=1\linewidth]{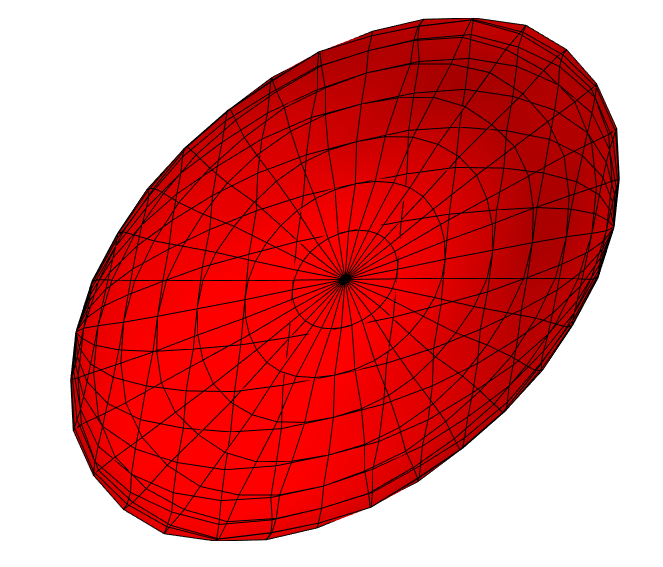}
\subcaption{$t=7.5$\label{fig:red_vol96_3}}
\end{subfigure}
\caption{Snapshots of capsule shape at different time instants  for shear flow simulation with $\Dot{\gamma}=1, E_{s}=2,E_{D}=200$. (a) $t=0$, (b) $t=0.25$ and (c) $t=7.5$. The shape shown in (c) is the terminal shape.  }
\label{fig:shear_redvol1}
\end{figure}

\begin{figure}[t!]
\begin{subfigure}{0.30\textwidth}
\includegraphics[width=1\linewidth]{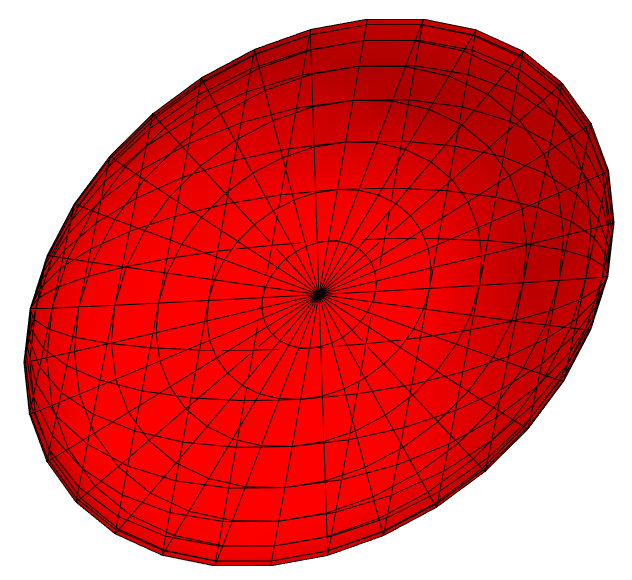}
\subcaption{$t=0.5$\label{fig:red_vol80_1}}
\end{subfigure}
\begin{subfigure}{0.23\textwidth}
\includegraphics[width=1\linewidth]{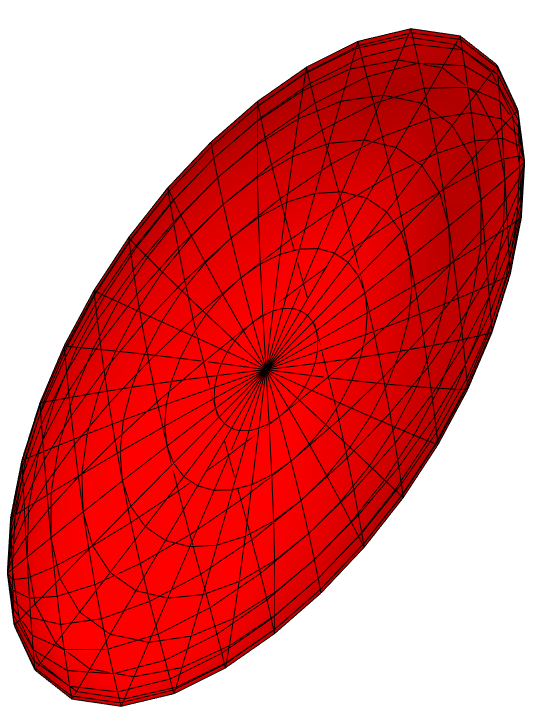}
\subcaption{$t=2$\label{fig:red_vol80_2}}
\end{subfigure}
\begin{subfigure}{0.35\textwidth}
\includegraphics[width=1\linewidth]{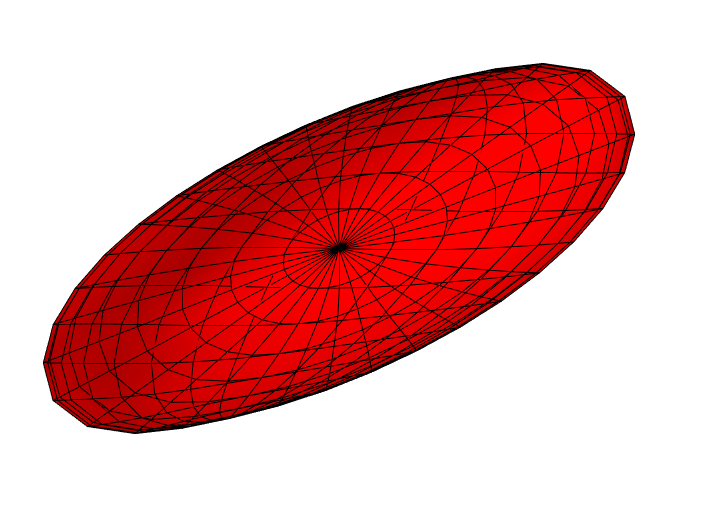}
\subcaption{$t=15$\label{fig:red_vol80_3}}
\end{subfigure}
\caption{Snapshots of capsule shape at different time instants  for shear flow simulation with $\Dot{\gamma}=0.5, E_{s}=2,E_{D}=1$. (a) $t=0.5$, (b) $t=2$ and (c) $t=15$. The shape shown in (c) is the terminal shape.  }
\label{fig:shear_redvol2}
\end{figure}

\begin{figure}[t!]
\begin{subfigure}{0.30\textwidth}
\includegraphics[width=1\linewidth]{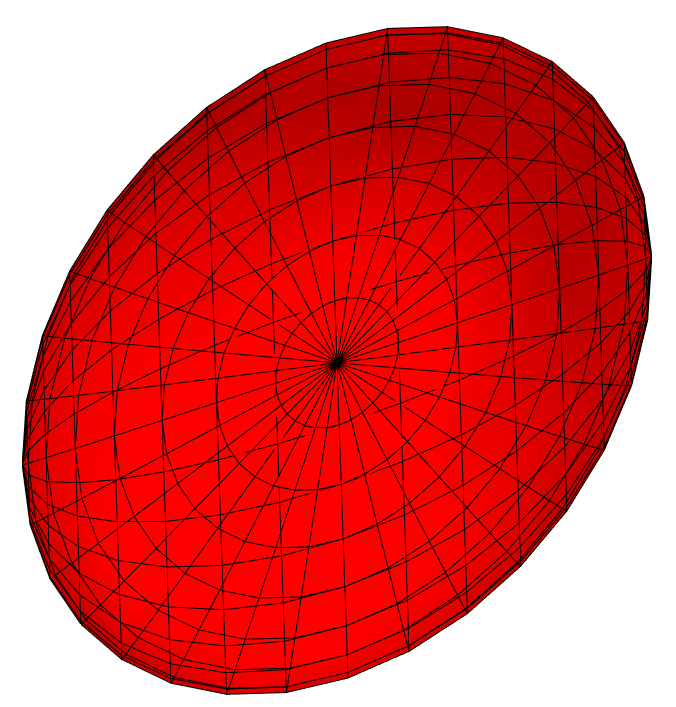}
\subcaption{$t=0.0417$\label{fig:red_vol60_1}}
\end{subfigure}
\begin{subfigure}{0.25\textwidth}
\includegraphics[width=1\linewidth]{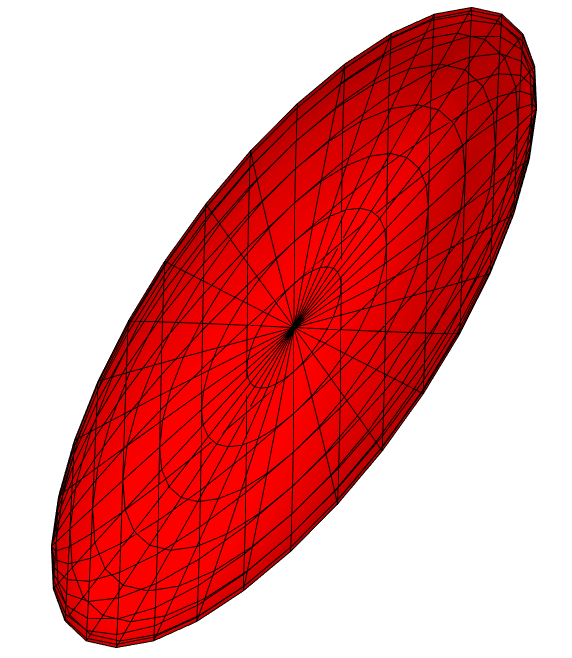}
\subcaption{$t=0.25$\label{fig:red_vol60_2}}
\end{subfigure}
\begin{subfigure}{0.35\textwidth}
\includegraphics[width=1\linewidth]{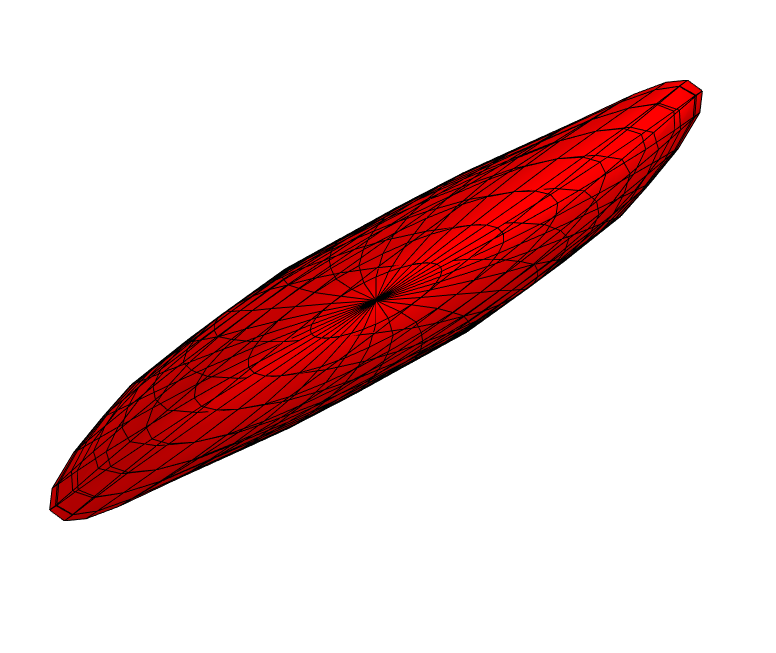}
\subcaption{$t=1.25$\label{fig:red_vol60_3}}
\end{subfigure}
\caption{Snapshots  of capsule shape at different time instants for shear flow simulation with $\Dot{\gamma}=6, E_{s}=2,E_{D}=1$. (a) $t=0.0417$, (b) $t=0.25$ and (c) $t=1.25$. The shape shown in (c) is the terminal shape.  }
\label{fig:shear_redvol3}
\end{figure}

Now, we present results for steady shapes of extensible capsule under background shear flow (see ~\cref{eq:shear_flow}). We start with a stress-free spherical capsule. Under linear shear flow, the capsule is known to take a terminal nearly-ellipsoidal shape and exhibit a stable tank treading motion~\cite{farutin2014, lac2004}. The terminal shape and inclination angle of its major axes with the flow direction depends on three key parameters, namely, the shear rate $\Dot{\gamma}$ of the flow, the shear modulus $E_{s}$ and the dilatation modulus $E_{D}$ of the membrane. We simulate these dynamics for different values of these parameters and plot the terminal inclination angles $\theta$ (with respect to the flow direction) in ~\cref{fig:shear_inclination}. We compare our results with the numerical results from ~\cite{farutin2014} and obtain good quantitative agreement. We also plot the evolution of Taylor asphericity $D_{a}$ with time in ~\cref{fig:shear_asphericity} and observe good agreement with the numerical results obtained using the spherical harmonics based scheme mentioned in ~\cite{veera2011}. ~\cref{fig:shear_redvol1,fig:shear_redvol2,fig:shear_redvol3} shows the terminal shapes of different reduced volumes obtained in shear flows using our code. Our code is able to resolve shapes for different reduced volumes obtained for a wide range of ratios of shear modulus $E_{s}$ and dilatation modulus $E_{D}$. 

\begin{figure}[t!]
\begin{subfigure}{0.55\textwidth}
\includegraphics[width=1\linewidth]{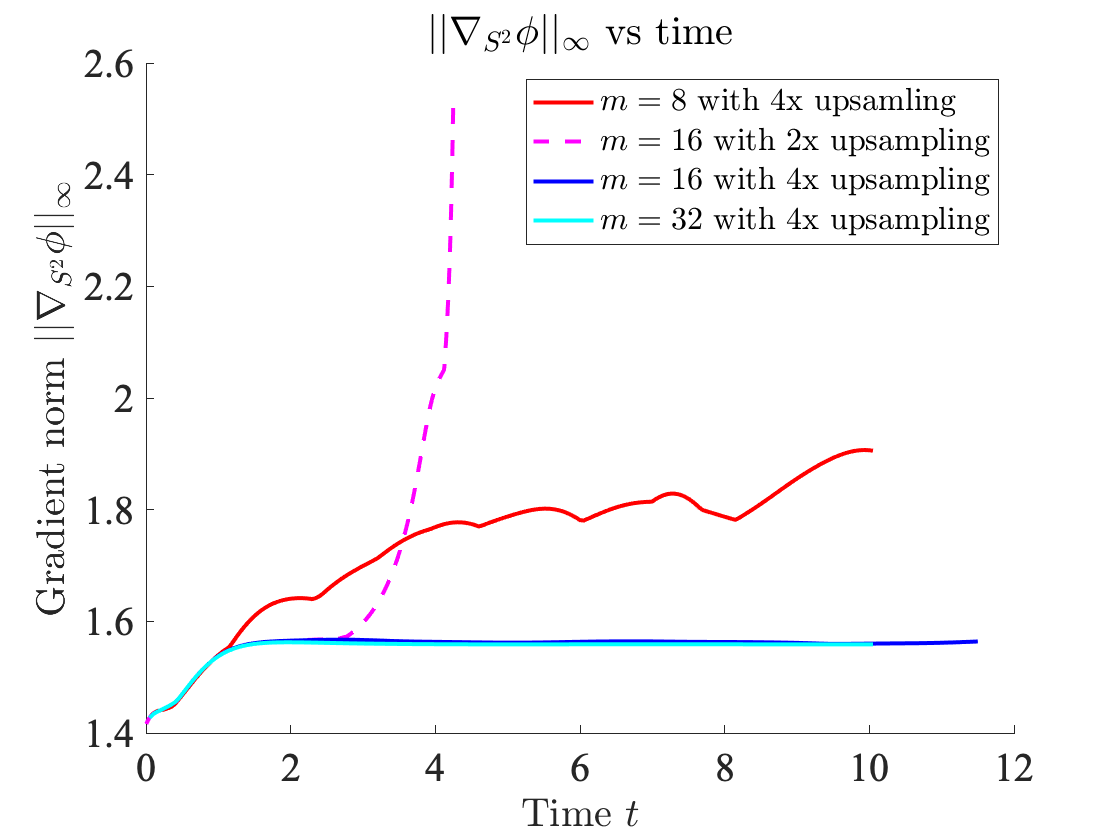}
\subcaption{Terminal reduced volume $\nu=0.85$. \label{fig:grad}}
\end{subfigure}%
\begin{subfigure}{0.55\textwidth}
\includegraphics[width=1\linewidth]{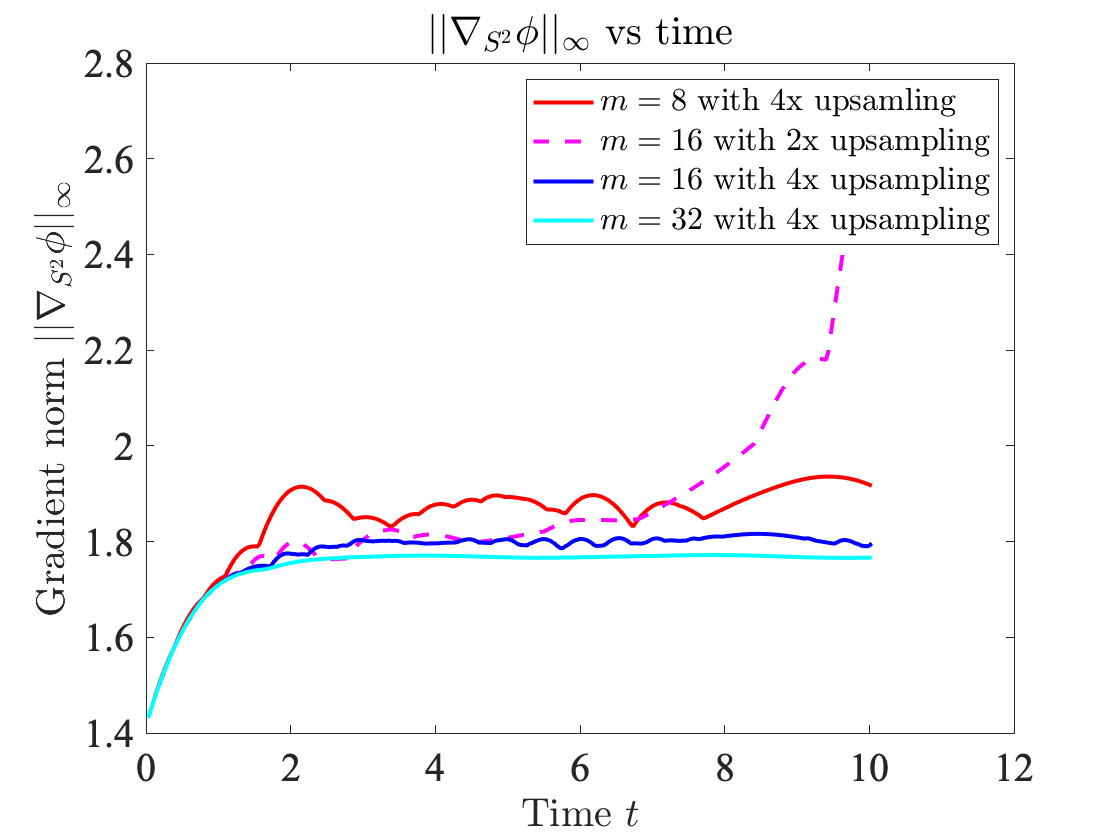}
\subcaption{Terminal reduced volume $\nu=0.65$. \label{fig:grad3}}
\end{subfigure}
\caption{ Plots of the $||\nabla_{\mathbb{S}^{2}} \phi||_{\infty}$ vs time $t$ in shear flow simulation where $\phi: \mathbb{S}^{2} \longrightarrow \gamma$ is the mapping from the unit sphere to the capsule $\gamma$. Initial shape is stress-free unit sphere with (a) $\Dot{\gamma}=1$, $E_{s}=2$ and $E_{D}=20$ and (b) $\Dot{\gamma}=1.5$, $E_{s}=2$ and $E_{D}=1$. The plots show the effectiveness of four times upsampling in doing stable long time horizon simulations while two times upsampling fails for $m=16$.    }
\label{fig:grads_shear}
\end{figure}

To further demonstrate the effectiveness of four times upsampling, we plot the norm of the gradient of mapping $\phi$ from a unit sphere to the capsule surface $\gamma$ with time in shear flow simulations for different grid orders $m$ in ~\cref{fig:grads_shear}. The gradient norm serves as a proxy for the stability of the capsule shape. As evident from the plots, we need high numerical accuracy to do long time horizon simulations. We observe that grid order $m=16$ with four times upsampling is sufficient to do shear flow simulations. Lower upsampling factors or lower grid order $m$ results in unstable gradients which blow up over long time scales as shown in the ~\cref{fig:grads_shear}. 

\subsection{Capsule in Poiseuille flow}
\begin{figure}[t!]
\begin{subfigure}{0.55\textwidth}
\includegraphics[width=1\linewidth]{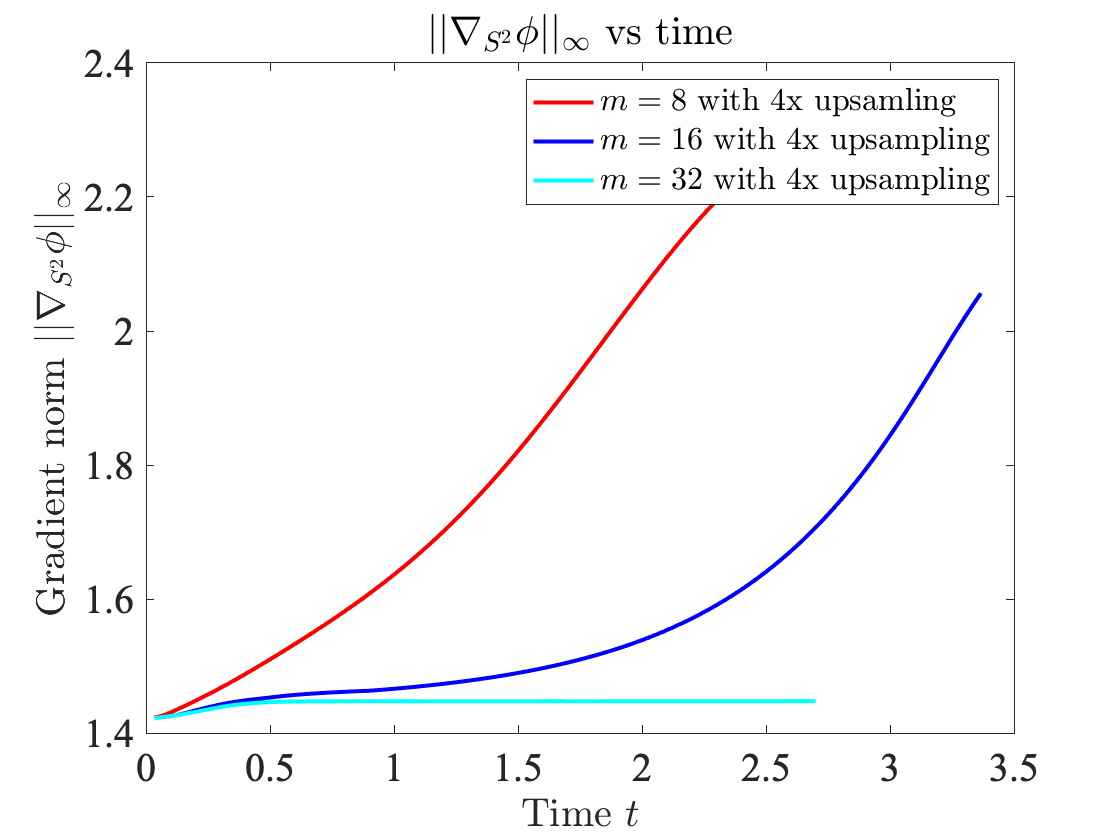}
\subcaption{Terminal reduced volume $\nu=0.96$. \label{fig:grad4}}
\end{subfigure}%
\begin{subfigure}{0.55\textwidth}
\includegraphics[width=1\linewidth]{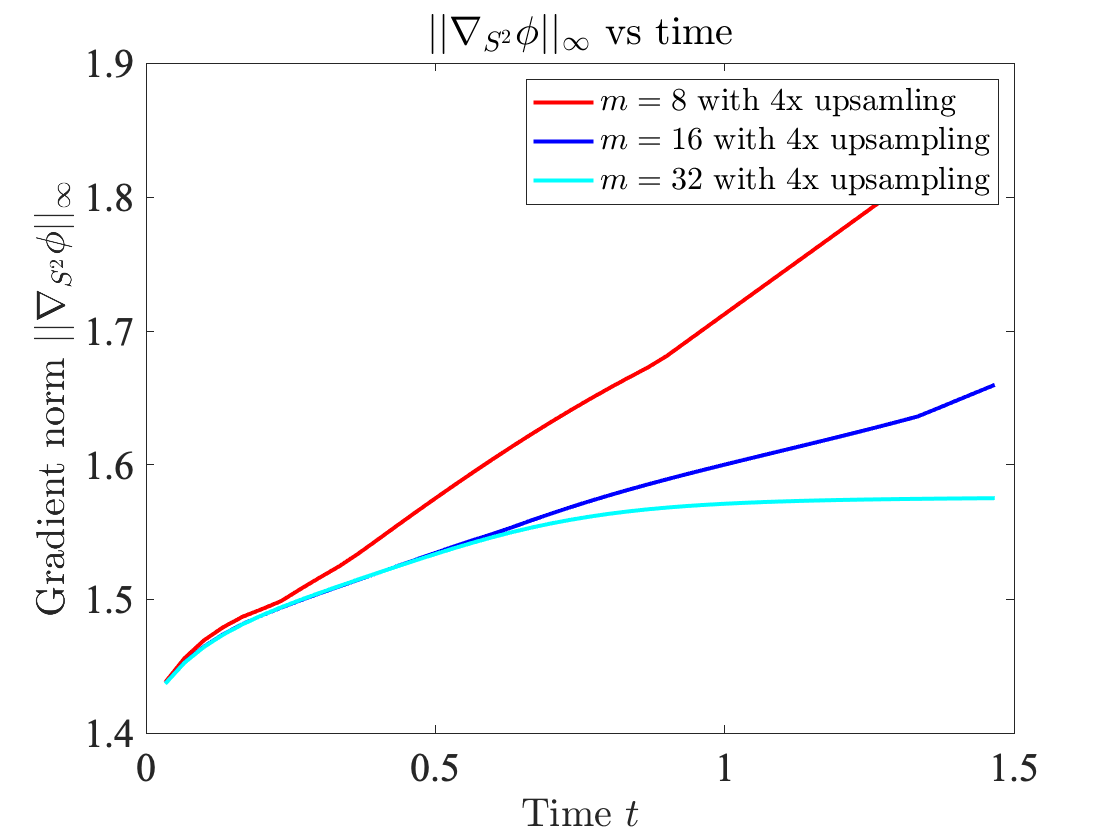}
\subcaption{Terminal reduced volume $\nu=0.86$. \label{fig:grad5}}
\end{subfigure}
\caption{ Plots of the $||\nabla_{\mathbb{S}^{2}} \phi||_{\infty}$ vs time $t$ in Poiseuille flow simulation where $\phi: \mathbb{S}^{2} \longrightarrow \gamma$ is the mapping from the unit sphere to the capsule $\gamma$. Initial shape is stress-free unit sphere with (a) $\alpha=1.5$, $R_{0}=5$, $E_{s}=2$ and $E_{D}=200$ and (b) $\alpha=1.5$, $R_{0}=5$, $E_{s}=2$ and $E_{D}=20$. Plots show that $m=32$ with four times upsampling gives stable simulations for Poiseuille flow while lower values of $m$ are not enough to resolves the shapes in Poiseuille flow.     }
\label{fig:grads_pois}
\end{figure}

\begin{figure}[t!]
\begin{subfigure}{0.30\textwidth}
\includegraphics[width=1\linewidth]{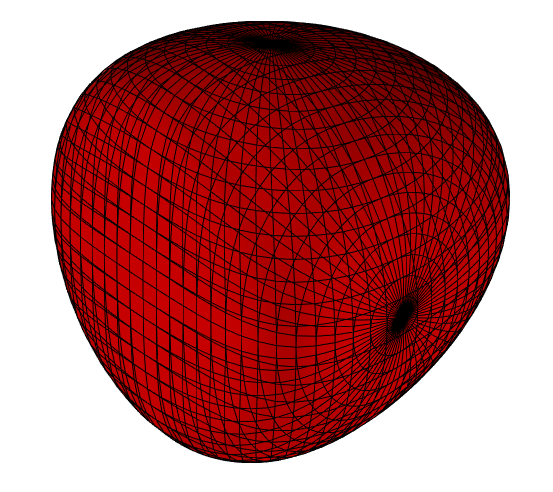}
\subcaption{$\nu=0.96$\label{fig:pois_red_vol96}}
\end{subfigure}
\begin{subfigure}{0.30\textwidth}
\includegraphics[width=1\linewidth]{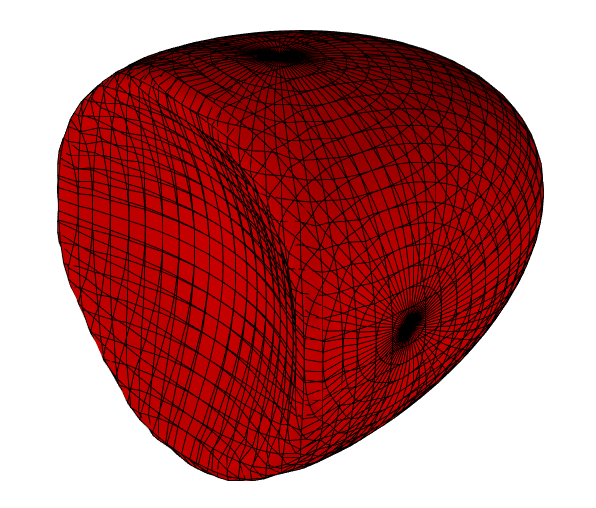}
\subcaption{$\nu=0.85$\label{fig:pois_red_vol80}}
\end{subfigure}
\begin{subfigure}{0.30\textwidth}
\includegraphics[width=1\linewidth]{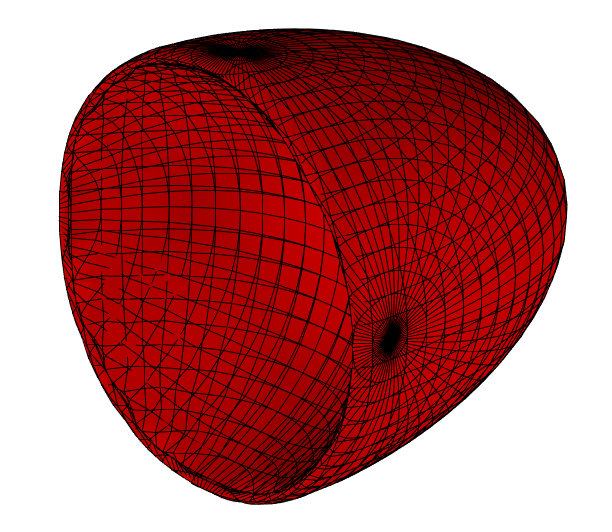}
\subcaption{$\nu=0.60$\label{fig:pois_red_vol65}}
\end{subfigure}
\caption{ Terminal parachute shapes of varying reduced volumes $\nu$ under Poiseuille flow. (a) Terminal shape of reduced volume $\nu=0.96$ obtained for Poiseuille flow simulation with $\alpha=1.5, E_{s}=2,E_{D}=200$. (b) Terminal shape of reduced volume  $\nu=0.85$ obtained for Poiseuille flow simulation with $\alpha=1.5, E_{s}=2,E_{D}=20$. (c) Terminal shape of reduced volume $\nu=0.6$ obtained for shear flow simulation with $\alpha=1.5, E_{s}=2,E_{D}=1$. We use grid order $m=32$ and the Poiseuille flow cross-section radius $R_{0}=5$ in these simulations . }
\label{fig:pois_redvol}
\end{figure}

In \cref{fig:pois_redvol}, we present the terminal shapes for a capsule under Poiseuille flow (see ~\cref{eq:poiseuille_flow}) for different membrane elasticity parameters leading to shapes of reduced volume as low as $\nu=0.6$. As for the shear flow simulations, we plot the norm of the gradient of mapping $\phi$ from a unit sphere to the capsule surface $\gamma$ with time for the Poiseuille flow simulations for different grid orders $m$ in ~\cref{fig:grads_pois}.  We observe that grid order $m=32$ with four times upsampling is sufficient to do Poiseuille flow simulations. Lower upsampling factors or lower grid order $m$ results in unstable gradients which blow up over long time scales as shown in the ~\cref{fig:grads_pois}.

\subsection{Timing results}
\begin{figure}[t!]
\centering
\begin{subfigure}{0.6\textwidth}
\includegraphics[width=1\linewidth]{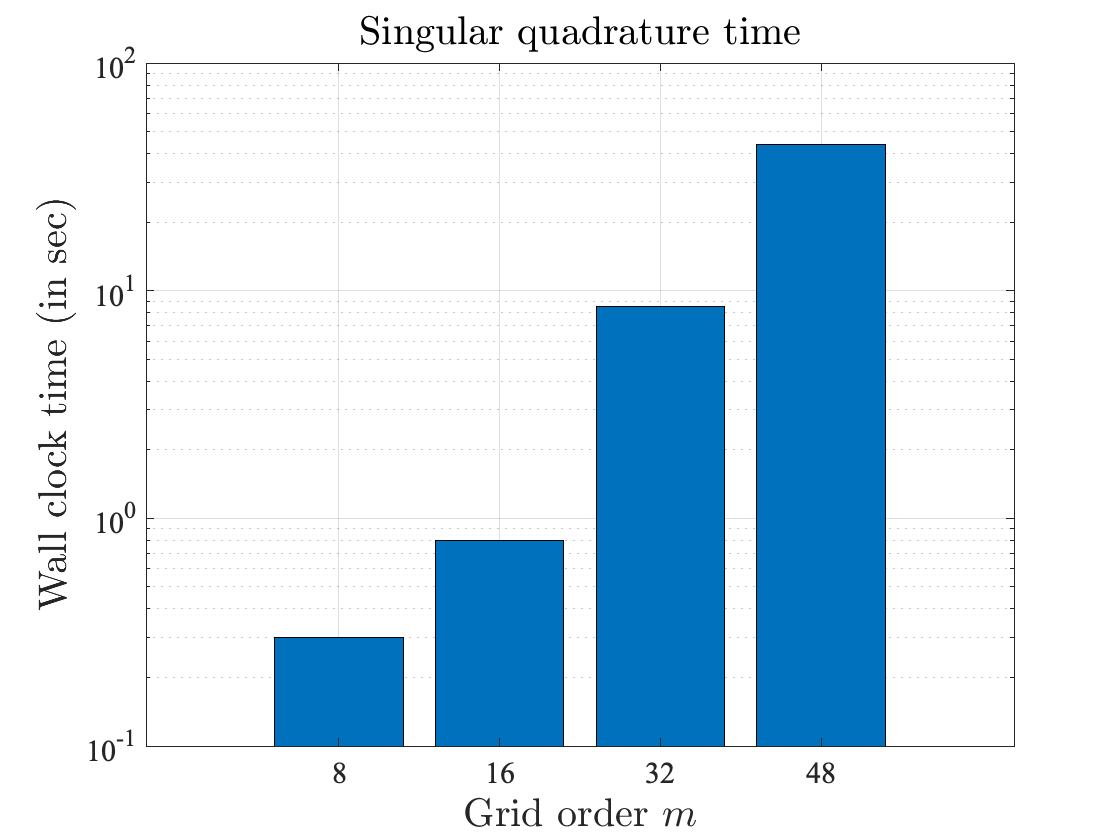}
\subcaption{\label{fig:SLtime}}
\end{subfigure}%
\begin{subfigure}{0.6\textwidth}
\includegraphics[width=1\linewidth]{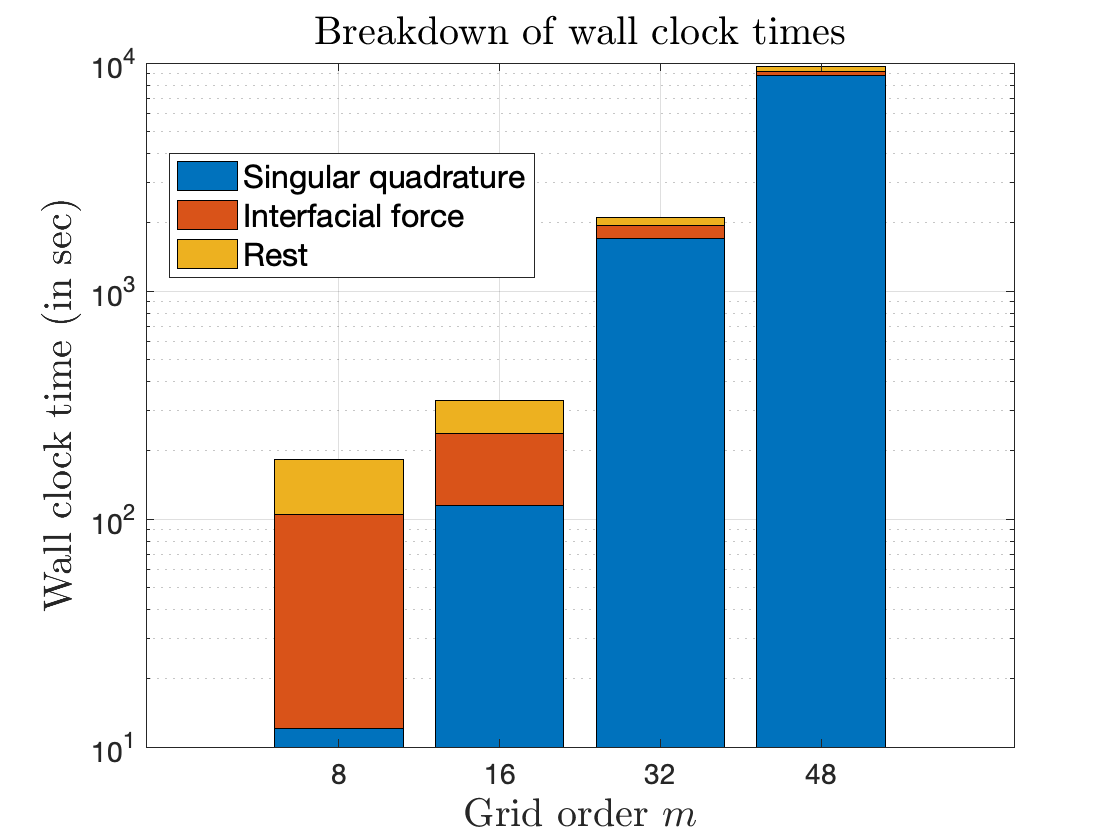}
\subcaption{\label{fig:time_dist}}
\end{subfigure}
\caption{ (a) GPU wall clock times for computing singular layer quadrature with four times upsampling  for different grid orders $m=8,16,32,48$. (b) Breakdown of wall clock times for shear flow simulation (with initial shape as unit sphere over time horizon $T = 0.1$) over different stages (computation of singular quadrature, computation of interfacial forces and the rest of the algorithm) for different grid orders $m$.   }
\label{fig:time}
\end{figure}

We compute the GPU wall clock times required for our scheme and give its breakdown over different stages of our scheme. To this end, we take an initially spherical capsule in stress free state (with $E_{s}=2, E_{D}=20$) and simulate its dynamics under shear flow with shear rate $\Dot{\gamma}=1$ till time $T=0.1$. We do these simulations for different grid order $m$ and plot the total time and its breakdown for different stages in ~\cref{fig:time_dist}. We observe that for $m\geq 16$, majority of the time is required to compute the singular quadrature. In fact, for $m \geq 32$ virtually all of the time is consumed in computing the quadrature. The wall clock times required to compute the singular quadrature once is also plotted separately in ~\cref{fig:SLtime}. We also provide a comparison of the wall clock time per time step taken by our GPU accelerated scheme with the spherical harmonics CPU code~\cite{veera2011} in ~\cref{app:wallclock}.  

\section{Fast multipole method based acceleration \label{sec:FMM}}
Our quadrature  does not involve a product quadrature and is amenable to  acceleration via fast multipole method (FMM) ~\cite{greengard2008}. A full FMM acceleration reduces the time complexity of our scheme to $O(N)$ allowing us to perform high resolution simulations. We do not implement a multilevel FMM for the purpose of this work and focus on showcasing the correctness and advantages of our numerical scheme using a simple single level FMM acceleration for our scheme. \da{Based on  the kernel independent FMM of~\cite{lexing2004}, we describe our scheme below. }

\begin{enumerate}
    \item We first use a $k$-means clustering algorithm to cluster together the discretization points on the capsule surface for a reasonable value of $k$. We will have $k$ clusters. We found that $k=100$ gives us relative accuracy of accuracy of $10^{-5}$, which is sufficient for the simulations in this paper.
    \item Each cluster is imagined to be enclosed by a surrounding cubical box, called an equivalent surface, with equivalent sources placed on a regular grid on the boundary of this cubical box. Thus, we have $k$ equivalent surfaces and let each equivalent surface have $N_{\mathrm{eq}}$ number of equivalent sources. We also consider a slightly bigger cube surrounding this enclosing cube which will be used as a check surface. The density on the equivalent sources lying on the equivalent surface is estimated by equating the potential on the check surface due to the actual point sources (\emph{i.e.,} discretization points) on the part of the capsule surface lying inside that box or equivalent surface. 
    \item We apply direct all pairs single layer potential calculation for the neighboring boxes while the equivalent sources are used to calculate the single layer potential for the discretization points in the far field (\emph{i.e.,} the discretization points inside non-neighboring boxes). 
\end{enumerate}

We present convergence and speedup results for calculating single layer on an ellipsoid of reduced volume $\nu=0.9$ in ~\cref{tab:ellipsoid_FMM}. While we do not see speedup for low values of $m$ because of FMM overhead computations, we see about two times speedup  for $m=96$ using our implementation of the single level FMM on GPU.

\begin{table}
\centering
\caption{Relative error and speedup for FMM accelerated simulation (with 4x upsampling) of a capsule suspended in Poiseuille flow as in ~\cref{fig:pois_red_vol65}. All pairs direct calculation for $m=32$ (with 4x upsampling) is used as the reference solution and the relative error is computed in the moment of inertia tensor of the capsule at the end of simulation. Time reported is wall clock time in seconds taken by one time step of the simulation. $t_{\mathrm{direct}}$ is direct calculation time and $t_{\mathrm{FMM}}$ is time taken with the FMM acceleration. Number of equivalent surfaces or boxes is fixed at 100. We vary the number of equivalent sources for an equivalent surface, denoted by $N_{eq}$, with $m$. }
\begin{tabular}{SSSSSSS} \toprule
    {$\textit{m}$} & {$N$} & {$N_{\mathrm{eq}}$} & {$t_{\mathrm{direct}}$}  & { $t_{\mathrm{FMM}}$} & {$\mathrm{Speedup}$} & $\epsilon_{\mathrm{FMM}}$\\ \midrule 
    32   & 5766  &  96  & 10.78  &    10.66  & 1.02  &  6E-3 \\
    64   & 23814  &  128  &  163.21  & 104.20 &  1.57  &  4E-5  \\ 
    96   & 54150  &  256  &   905.31     &  501.76 &  1.81  &  5E-6 \\
    \bottomrule

\end{tabular}
\subcaption{ \label{tab:ellipsoid_FMM}}.  
\end{table}

\section{Conclusions \label{sec:conclusion}}
In this work, we described a novel numerical scheme to simulate Stokesian particulate flows. We described an overlapping patch based discretization of the surfaces diffeomorphic to the unit sphere. Our numerical scheme uses the regularized Stokes kernels ~\cite{beale2019} and finite differences on overset grid to calculate the Stokes layer potential and interfacial elastic forces. We presented a battery of results for the verification of our numerical scheme using results from previous literature. We used our numerical scheme to simulate an extensible capsule in suspended in Stokes flow and validated our simulations. Our numerical scheme is a fourth order convergent scheme that is $O(N^{2})$ in work complexity where $N$ is the number of discretization points and can be accelerated to run with $O(N)$ work complexity using a multilevel FMM. This is much better than the asymptotic work complexity of the spherical harmonics based spectral scheme~\cite{veera2011}. Our scheme also allows for independent control over local resolution due to patch based parameterization of surface. We used GPU acceleration to demonstrate the ability of our code to simulate the complex shapes with high resolution. 

In this paper, we implemented a single level FMM to demonstrate FMM based speedup. A GPU based implementation of a multilevel FMM will further enable us to do highly accurate simulations in a reasonable time and help us better study the Stokesian particulate flows. We leave it as the subject of our future work.

\section{Funding}
Dhwanit Agarwal and George Biros' research was supported in part by UT-Austin CoLab and UTEN
Partnership, and the Portuguese Science and Technology Foundation under funding
~\#201801976/UTA18-001217. 

\section{Conflict of interest}
The authors do not have any conflicts of interest.

\appendix
\section{Appendix }
\subsection{Transition maps \label{app:transition}}
Here, we list the expressions for the transition maps $\tau_{ij}$ for $i=1,..,6$, for the parameterization we use. The expressions are given as follows:
\begin{align}
    \tau_{12}(u,v) &= (u,\pi +v), \tau_{13}(u,v) = (u,\frac{\pi}{2}+v), \tau_{14}(u,v) = (u,v-\frac{\pi}{2}),\\
    \tau_{15}(u,v) &= \left(\cosinv(-\sin{u}\sin{v}),\cosinv \left(\frac{\sin{u}\cos{v}}{\sqrt{\sinsq{u} \cossq{v} + \cossq{u}}}\right) \right ),\\
    \tau_{16}(u,v) &= \left(\cosinv(\sin{u}\sin{v}),\cosinv \left(\frac{\sin{u}\cos{v}}{\sqrt{\sinsq{u} \cossq{v} + \cossq{u}}}\right) \right ).
\end{align}

The transition maps for the next three patches are given below: 

\begin{align}
    \tau_{21}(u,v) &= (u,\pi +v), \tau_{23}(u,v) = (u,v-\frac{\pi}{2}), \tau_{24}(u,v) = (u,v+\frac{\pi}{2}),\\
    \tau_{25}(u,v) &= \left(\cosinv(\sin{u}\sin{v}),\cosinv \left(\frac{-\sin{u}\cos{v}}{\sqrt{\sinsq{u} \cossq{v} + \cossq{u}}}\right) \right ),\\
    \tau_{26}(u,v) &= \left(\cosinv(-\sin{u}\sin{v}),\cosinv \left(\frac{-\sin{u}\cos{v}}{\sqrt{\sinsq{u} \cossq{v} + \cossq{u}}}\right) \right ),\\
    \tau_{31}(u,v) &= (u,v- \pi/2), \tau_{32}(u,v) = (u,v+\frac{\pi}{2}), \tau_{34}(u,v) = (u,v+\pi),\\
    \tau_{35}(u,v) &= \left(\cosinv(\sin{u}\cos{v}),\cosinv \left(\frac{\sin{u}\sin{v}}{\sqrt{\sinsq{u} \sinsq{v} + \cossq{u}}}\right) \right ),\\
    \tau_{36}(u,v) &= \left(\cosinv(-\sin{u}\cos{v}),\cosinv \left(\frac{\sin{u}\sin{v}}{\sqrt{\sinsq{u} \sinsq{v} + \cossq{u}}}\right) \right ),\\
    \tau_{41}(u,v) &= (u,v+\pi/2), \tau_{42}(u,v) = (u,v-\frac{\pi}{2}), \tau_{43}(u,v) = (u,v+\pi),\\
    \tau_{45}(u,v) &= \left(\cosinv(-\sin{u}\cos{v}),\cosinv \left(\frac{-\sin{u}\sin{v}}{\sqrt{\sinsq{u} \sinsq{v} + \cossq{u}}}\right) \right ),\\
    \tau_{46}(u,v) &= \left(\cosinv(\sin{u}\cos{v}),\cosinv \left(\frac{-\sin{u}\sin{v}}{\sqrt{\sinsq{u} \sinsq{v} + \cossq{u}}}\right) \right ).
\end{align}
The transition maps for $\mathcal{P}_{5}^{0}$ and $\mathcal{P}_{6}^{0}$ are given as follows:

\begin{align}
    \tau_{51}(u,v) &= \left(\cosinv(\sin{u}\sin{v}),\cosinv \left(\frac{\sin{u}\cos{v}}{\sqrt{\sinsq{u} \cossq{v} + \cossq{u}}}\right) \right ),\\
    \tau_{52}(u,v) &= \left(\cosinv(\sin{u}\sin{v}),\cosinv \left(\frac{-\sin{u}\cos{v}}{\sqrt{\sinsq{u} \cossq{v} + \cossq{u}}}\right) \right ),\\
    \tau_{53}(u,v) &= \left(\cosinv(\sin{u}\sin{v}),\cosinv \left(\frac{\cos{u}}{\sqrt{\sinsq{u} \cossq{v} + \cossq{u}}}\right) \right ),\\
    \tau_{54}(u,v) &= \left(\cosinv(\sin{u}\sin{v}),\cosinv \left(\frac{-\cos{u}}{\sqrt{\sinsq{u} \cossq{v} + \cossq{u}}}\right) \right ),\\
    \tau_{54}(u,v) &= (\pi-u,v),\\
    \tau_{61}(u,v) &= \left(\cosinv(-\sin{u}\sin{v}),\cosinv \left(\frac{\sin{u}\cos{v}}{\sqrt{\sinsq{u} \cossq{v} + \cossq{u}}}\right) \right ),\\
    \tau_{62}(u,v) &= \left(\cosinv(-\sin{u}\sin{v}),\cosinv \left(\frac{-\sin{u}\cos{v}}{\sqrt{\sinsq{u} \cossq{v} + \cossq{u}}}\right) \right ),\\
    \tau_{63}(u,v) &= \left(\cosinv(-\sin{u}\sin{v}),\cosinv \left(\frac{\cos{u}}{\sqrt{\sinsq{u} \cossq{v} + \cossq{u}}}\right) \right ),\\
    \tau_{64}(u,v) &= \left(\cosinv(-\sin{u}\sin{v}),\cosinv \left(\frac{-\cos{u}}{\sqrt{\sinsq{u} \cossq{v} + \cossq{u}}}\right) \right ),\\
    \tau_{65}(u,v) &= (\pi-u,v).
\end{align}

\subsection{Surface derivative formulas \label{app:formulas}}
Here, we list the formulas for the first fundamental form coefficients $E,F,G$, the unit normal $\bn$, the second fundamental form coefficients $L,M,N$, the mean curvature $H$, the Gaussian curvature $K$, the surface divergence (of a vector field $\bg$) and the surface gradient (of a scalar function $g$) in terms of the local parameterization $\bx(u,v):D \subset \mathbb{R}^{2} \longrightarrow  \gamma$. We use these quantities in the computation of interfacial force and the verification of surface derivatives.  

\begin{table}[h!]
\centering
\caption{Formulas used for computing various surface derivative quantities. \label{tab:diff_formulas}}
 \begin{tabular}{|p{2cm}  p{4cm} p{2cm} p{4cm}|}
 \hline
 Symbol & Definition & Symbol & Definition \\ [0.5ex] 
 \hline\hline
  $E$ & $\bx_{u} \cdot \bx_{u}$ &    $M$   & $\bx_{uv}\cdot \bn$  \\
  \hline
  $F$ & $\bx_{u} \cdot \bx_{v}$ &    $N$   & $\bx_{vv}\cdot \bn$   \\
  \hline
  $G$ &  $\bx_{v} \cdot \bx_{v}$&     $H$ & $\frac{EN-2FM+GL}{2W^{2}}$ \\
  \hline
  $W$ & $\sqrt{EG-F^{2}}$ &      $\nabla_{\gamma} g$   &   $\frac{G\bx_{u}-F\bx_{v}}{W^{2}}g_{u} + \frac{E\bx_{v}-F\bx_{u}}{W^{2}}g_{v}$    \\
  \hline
  $\bn$ & $ \frac{\bx_{u}\times \bx_{v}}{W}$  & $\nabla_{\gamma} \cdot \bg$  & $\frac{G\bg_{u}-F\bg_{v}}{W^{2}}\bx_{u} + \frac{E\bg_{v}-F\bg_{u}}{W^{2}}\bx_{v}$    \\ 
  \hline
  $L$ & $\bx_{uu}\cdot \bn$ &  $K$ & $\frac{LN - M^{2}}{W^{2}}$     \\
\hline 
  \end{tabular}
\end{table}

\subsection{The choice of parameter $r_{0}$ for the partition of unity \label{app:err_pou_overlap}}
Here, we briefly discuss our choice of parameter $r_{0}$ in the construction of partition of unity in ~\cref{sub:surfparam}. We note here that for the specific parameterization of the unit sphere we use (see ~\cref{eq:param_maps}), we require   $r_{0}>\frac{3\pi}{12}$ to ensure that every point on $\mathbb{S}^{2}$ belongs to the support of $\psi_{i}^{0}$ for at least one $i\in \{1,\ldots,6\}$. If this is not the case, then $\{\psi_{i}^{0}\}_{i=1}^{6}$ ceases to be a partition of unity on the unit sphere. Additionally, we need $r_{0}<\frac{\pi}{2}$ so that the support of each $\psi_{i}^{0}$ is compactly contained in $\mathcal{P}_{i}^{0}$. We tabulate some numerical results ~\cref{tab:sph_FDM} for the relative errors in the surface derivatives on the unit sphere for different values of $\pi > r_{0} > \frac{3\pi}{12}$. We find that $ r_{0}=\frac{5\pi}{12}$ gives optimal accuracy among the values of $r_{0}$ we experimented with. 

\begin{table}
\centering
\caption{Relative errors in surface derivatives. We compare the relative errors for different values of $r_{0}$. }
\begin{tabular}{SSSS} \toprule
    {$m$}  & {$r_{0}/\pi$} & {$\epsilon_{\bn}$} & {$\epsilon_{H}$}    \\ \midrule
    8   &  5.5/12    &  2E-4  &   1.3E-3     \\
    8   &  5/12    &   2E-4   &  1.5E-3      \\
    8   &  4/12     &  4E-4   &  1.6E-3     \\ \midrule
    16  & 5.5/12     &  8E-6   &   1.7E-5     \\
    16  &  5/12     &  7E-6    &   1.8E-5     \\ 
    16  &  4/12     &   7E-6   &    1.9E-5  \\ \midrule
    32 &  5.5/12     &  4E-7   &   1.3E-6     \\ 
    32  &  5/12     &  4E-7    &   1.7E-6     \\ 
    32  &  4/12     &   4E-7   &    1.9E-6   \\ \bottomrule

\end{tabular}
 \label{tab:sph_FDM}
\end{table}

\subsection{Effect of the blending process in calculating derivatives \label{app:blending_evidence}}
Here, we provide the relative errors in the computing normals $n$ and the mean curvature $H$ for the unit sphere with and without blending process discussed in ~\cref{sub:surfderiv}. The results are tabulated in ~\cref{tab:blending_evidence}. The results show that blending improves the accuracy of derivatives.  

\begin{table}
\centering
\caption{ Relative error in computing surface normals $n$ and the mean curvature $H$ for unit sphere with and without the blending. $\epsilon^{\mathrm{b}}$ are the errors with blending and $\epsilon^{\mathrm{nb}}$  are the errors without the blending process in computing derivatives. \label{tab:blending_evidence}}
\begin{tabular}{SSSSS} \toprule
    {$m$}  & {$\epsilon_{\bn}^{\mathrm{b}}$} & {$\epsilon_{H}^{\mathrm{b}}$} & {$\epsilon_{\bn}^{\mathrm{nb}}$} & {$\epsilon_{H}^{\mathrm{nb}}$}\\ \midrule
    8  & 2E-4 & 1.3E-3 & 3E-3 & 2E-2    \\
    16 & 8E-6 & 1.7E-5 & 9E-5 & 2E-3     \\    
    32 & 4E-7 & 1.3E-6 & 3E-6 & 3E-4     \\  \bottomrule

\end{tabular}
\end{table}

\subsection{Surface derivative and singular integration errors for shear flow and Poiseuille flow terminal shapes \label{app:err_shear_pois_shapes}}
To further verify the convergence and accuracy of our derivative and integration schemes, we upsample the low reduced volume shapes obtained in ~\cref{fig:red_vol60_3} and ~\cref{fig:pois_red_vol65} to $p=64$ spherical harmonics and use the spherical harmonics derivatives and Graham-Sloan quadrature ~\cite{veera2011} as the reference values to study the convergence of our schemes. We report the relative errors for these shapes in ~\cref{tab:err_shear_pois}.

\begin{table}
\centering
\caption{Relative errors in the derivative scheme and singular quadrature for the shapes in ~\cref{fig:red_vol60_3} and ~\cref{fig:pois_red_vol65} using our scheme (at the top for different grid orders $m$ with four times upsampling to compute quadrature ).  The shapes are upsampled to spherical harmonics degree $p=64$ and the error is computed relative to those values. Surface divergence is computed for the smooth function $\bbf(x,y,z)=(x^{2},y^{2},z^{2})$ on the surface. \label{tab:err_shear_pois}  }
\begin{tabular}{SSSSSSS} \toprule
    {$\textit{m}$} & {$N$} &  {$\epsilon^{up}_{\mathcal{S}[\bbf]}$} & {$\epsilon_{\bn}$}  & { $\epsilon_{H}$} & {$\epsilon_{K}$} & { $\epsilon_{div_{\gamma}}$}\\ \midrule
    8   & 294  & 5E-2 &    9E-2  & 5E-1 &  1E0    & 9E-2      \\
    16  & 1350 & 2E-2  &   6E-2   &  2E-1 &  2E-1  & 5E-2    \\    
    32  & 5766 & 5E-3  &    9E-3  & 8E-2  &  9E-2  & 9E-3     \\
    64  & 23814 & 5E-4  &    9E-4  &  8E-3  &  1E-2 &  9E-4 \\  \bottomrule
\end{tabular}
\subcaption{Relative error in singular quadrature and derivatives for the shape in ~\cref{fig:red_vol60_3} of reduced volume $\nu=0.6$.  }
\vspace{1cm}
\begin{tabular}{SSSSSSS} \toprule
    {$\textit{m}$} & {$N$} &  {$\epsilon^{up}_{\mathcal{S}[\bbf]}$} & {$\epsilon_{\bn}$}  & { $\epsilon_{H}$} & {$\epsilon_{K}$} & { $\epsilon_{div_{\gamma}}$}\\ \midrule
    8   & 294  & 8E-2 &    9E-2  & 7E-1 &  8E-1    & 9E-2      \\
    16  & 1350 & 2E-2  &   5E-2   &  2E-1 &  3E-1  & 6E-2    \\    
    32  & 5766 & 3E-3  &    8E-3  & 6E-2  &  7E-2  & 9E-3     \\
    64  & 23814 & 6E-4  &    9E-4   &  1E-2  &  1E-2 &  8E-4 \\  \bottomrule

\end{tabular}
\subcaption{Relative error in singular quadrature and derivatives for the shape ~\cref{fig:pois_red_vol65} of reduced volume $\nu=0.60$.   }
\end{table}

\subsection{Numerical errors for different reduced volume $\nu$ \label{app:err_redvol}}
Here, we tabulate errors in computing the singular quadrature and surface derivatives using our numerical scheme. We report these errors for ellipsoids of varying reduced volumes $\nu \in \{0.92,0.78,0.5,0.3\}$ in ~\cref{tab:err_redvol} and ~\cref{tab:err_redvol2}. The tables demonstrate that the accuracy deteriorates with decreasing reduced volumes. In general, it is harder to do long time horizon simulations of capsules with lower reduced volumes due to this reason. 
\begin{table}
\centering
\caption{Relative errors in the derivative scheme and singular quadrature for ellipsoids of reduced volume (a) $\nu=0.92$ and (b) $\nu=0.78$  using our scheme (at the top for different grid orders $m$ with four times upsampling to compute quadrature ) and using spherical harmonics scheme~\cite{veera2011} (at the bottom for different spherical harmonics degrees $p$).  Spherical harmonics results for $p=64$ are used as true values and the error is computed relative to those values. Surface divergence is computed for the smooth function $\bbf(x,y,z)=(x^{2},y^{2},z^{2})$ on the surface. \label{tab:err_redvol}  }
\begin{tabular}{SSSSSSS} \toprule
    {$\textit{m}$} & {$N$} &  {$\epsilon^{up}_{\mathcal{S}[\bbf]}$} & {$\epsilon_{\bn}$}  & { $\epsilon_{H}$} & {$\epsilon_{K}$} & { $\epsilon_{div_{\gamma}}$}\\ \midrule
    8   & 294  & 7E-3 &    5E-3  & 4E-3 &  6E-3    & 4E-2      \\
    16  & 1350 & 4E-4  &   2E-4   &  5E-4 &  1E-3  & 3E-3    \\    
    32  & 5766 & 2E-5  &    1E-5  & 3E-5  &  6E-5  & 1E-4     \\
    64  & 23814 & 1E-6  &    7E-7  &  2E-6  &  3E-6 &  1E-5 \\  \midrule \midrule
    {$\textit{p}$} & {$N$}  &  {$\epsilon_{\mathcal{S}[\bbf]}$} & {$\epsilon_{\bn}$}  & { $\epsilon_{H}$} & {$\epsilon_{K}$} & { $\epsilon_{div_{\gamma}}$}\\ \midrule
    8    &  144 & 4E-5 &    3E-3  & 2E-3 &  6E-3  &   6E-3  \\
    16   & 544 & 3E-8  &   1E-5   &  1E-5 &  4E-5 &   2E-5 \\    
    24  & 1200 & 5E-11  &    4E-8  & 6E-8  &  2E-7 &   1E-7 \\
    32   & 2112 & 6E-13  &    5E-11   & 8E-10  & 4E-9 & 3E-9 \\  \midrule \midrule    
\end{tabular}
\subcaption{Relative error in singular quadrature and derivatives for an ellipsoid of reduced volume $\nu=0.9$, given by $\frac{x^{2}}{a^{2}} + \frac{y^{2}}{b^{2}} + \frac{z^{2}}{c^{2}} = 1$, where $a=0.6, b=1, c=1$.  }
\vspace{1cm}
\begin{tabular}{SSSSSSS} \toprule
    {$\textit{m}$} & {$N$} &  {$\epsilon^{up}_{\mathcal{S}[\bbf]}$} & {$\epsilon_{\bn}$}  & { $\epsilon_{H}$} & {$\epsilon_{K}$} & { $\epsilon_{div_{\gamma}}$}\\ \midrule
    8   & 294  & 7E-3 &    3E-2  & 2E-2 &  4E-2    & 4E-2      \\
    16  & 1350 & 4E-4  &   9E-4   &  3E-3 &  5E-3  & 9E-4    \\    
    32  & 5766 & 2E-5  &    8E-5  & 3E-4  &  5E-4  & 1E-4     \\
    64  & 23814 & 1E-6  &    4E-6   &  1E-5  &  3E-5 &  9E-6 \\  \midrule \midrule
    {$\textit{p}$} & {$N$}  &  {$\epsilon_{\mathcal{S}[\bbf]}$} & {$\epsilon_{\bn}$}  & { $\epsilon_{H}$} & {$\epsilon_{K}$} & { $\epsilon_{div_{\gamma}}$}\\ \midrule
    8    &  144 & 7E-4 &    2E-2  & 3E-2 &  6E-2  &   4E-2  \\
    16   & 544 & 2E-6  &   6E-4   &  1E-3 &  3E-3 &   1E-3 \\    
    24  & 1200 & 3E-8  &    2E-5  & 5E-5  &  1E-4 &   5E-5 \\
    32   & 2112 & 5E-10  &    1E-7   &  6E-7  &  9E-6 & 4E-6 \\  \midrule \midrule    
\end{tabular}
\subcaption{Relative error in singular quadrature and derivatives for an ellipsoid of reduced volume $\nu=0.78$, given by $\frac{x^{2}}{a^{2}} + \frac{y^{2}}{b^{2}} + \frac{z^{2}}{c^{2}} = 1$, where $a=0.4, b=1, c=1$.  }

\end{table}

\begin{table}
\centering
\caption{Relative errors in the derivative scheme and singular quadrature for ellipsoids of reduced volume (a) $\nu=0.5$ and (b) $\nu=0.3$  using our scheme (at the top for different grid orders $m$ with four times upsampling to compute quadrature ) and using spherical harmonics scheme~\cite{veera2011} (at the bottom for different spherical harmonics degrees $p$).  Spherical harmonics results for $p=64$ are used as true values and the error is computed relative to those values. Surface divergence is computed for the smooth function $\bbf(x,y,z)=(x^{2},y^{2},z^{2})$ on the surface. \label{tab:err_redvol2}  }
\begin{tabular}{SSSSSSS} \toprule
    {$\textit{m}$} & {$N$} &  {$\epsilon^{up}_{\mathcal{S}[\bbf]}$} & {$\epsilon_{\bn}$}  & { $\epsilon_{H}$} & {$\epsilon_{K}$} & { $\epsilon_{div_{\gamma}}$}\\ \midrule
    8   & 294  & 2E-2 &    2E-1  & 2E-1 &  3E-1    & 4E-2      \\
    16  & 1350 & 4E-3  &   9E-3   &  3E-2 &  6E-2  & 3E-3    \\    
    32  & 5766 & 7E-4  &    8E-4  & 4E-3  &  6E-3  & 1E-4     \\
    64  & 23814 & 6E-5  &    9E-5  &  4E-4  &  7E-4 &  1E-5 \\  \midrule \midrule
    {$\textit{p}$} & {$N$}  &  {$\epsilon_{\mathcal{S}[\bbf]}$} & {$\epsilon_{\bn}$}  & { $\epsilon_{H}$} & {$\epsilon_{K}$} & { $\epsilon_{div_{\gamma}}$}\\ \midrule
    8    &  144 & 9E-3 &    9E-2  & 3E-1 &  4E-1  &   2E-3 \\
    16   & 544 & 2E-4  &   1E-2   &  6E-2 &  9E-2 &   4E-2 \\    
    24  & 1200 & 1E-5  &    3E-3  & 1E-2  &  2E-2 &   7E-3 \\
    32   & 2112 & 2E-6  &    1E-4   &  1E-3  &  5E-3 & 6E-4 \\  \midrule \midrule    
\end{tabular}
\subcaption{Relative error in singular quadrature and derivatives for an ellipsoid of reduced volume $\nu=0.5$, given by $\frac{x^{2}}{a^{2}} + \frac{y^{2}}{b^{2}} + \frac{z^{2}}{c^{2}} = 1$, where $a=0.2, b=1, c=1$.  }
\vspace{1cm}
\begin{tabular}{SSSSSSS} \toprule
    {$\textit{m}$} & {$N$} &  {$\epsilon^{up}_{\mathcal{S}[\bbf]}$} & {$\epsilon_{\bn}$}  & { $\epsilon_{H}$} & {$\epsilon_{K}$} & { $\epsilon_{div_{\gamma}}$}\\ \midrule
    8   & 294  & 4E-2 &    4E-1  & 5E-1 &  6E-1    & 4E-2      \\
    16  & 1350 & 1E-2  &   9E-2   &  2E-1 &  3E-1  & 3E-3    \\    
    32  & 5766 & 4E-3  &    2E-2  & 4E-2  &  7E-2  & 2E-3     \\
    64  & 23814 & 1E-3  &    2E-3  &  8E-3  &  9E-3 &  9E-4 \\  \midrule \midrule
    {$\textit{p}$} & {$N$}  &  {$\epsilon_{\mathcal{S}[\bbf]}$} & {$\epsilon_{\bn}$}  & { $\epsilon_{H}$} & {$\epsilon_{K}$} & { $\epsilon_{div_{\gamma}}$}\\ \midrule
    8    &  144 & 4E-2 &    2E-1  & 1E0 &  1E0  &   4E-1  \\
    16   & 544 & 3E-3  &   1E-1   &  3E-1 &  4E-1 &   2E-1 \\    
    24  & 1200 & 5E-4  &    4E-2  & 1E-1  &  1E-1 &   6E-2 \\
    32   & 2112 & 4E-5  &    5E-3   &  3E-2  & 4E-2 & 5E-3 \\  \midrule \midrule    
\end{tabular}
\subcaption{Relative error in singular quadrature and derivatives for an ellipsoid of reduced volume $\nu=0.3$ given by $\frac{x^{2}}{a^{2}} + \frac{y^{2}}{b^{2}} + \frac{z^{2}}{c^{2}} = 1$, where $a=0.1, b=1, c=1$.  }

\end{table}

\subsection{Sensitivity of the singular quadrature scheme to the regularization parameter $\delta$  \label{app:sensitivity}}
Here, we discuss the effect of choice of regularization parameter $\delta$ on the singular quadrature accuracy. We experiment with different values of the constant in $C$ in ~\cref{eq:regparam}. We tabulate the error in singular quadrature  for different values of $C$ in ~\cref{tab:sing_int_delta}. We find out that $C=1$ works best for the parameterization we have chosen. We also tabulate in the same table the results for fixed global regularization parameter $\delta = C h_{m}$ and show that our patch-dependent regularization parameter gives better accuracy and convergence compared to the fixed regularization parameter. 
\begin{table}
\centering
\caption{Sensitivity of singular quadrature with respect to the regularization parameter $\delta$. First three columns show the relative errors for the patch-dependent $\delta = C\delta^{*}$ where $C$ is a constant. The last three columns show the relative errors for fixed global regularization parameter $\delta = C h_{m}$.   \label{tab:sing_int_delta}}
\begin{tabular}{SSSSSSS} \toprule
    {$m$}  & {$\delta=0.5\delta^{*}$  } &  {$\delta=\delta^{*}$} & {$\delta=2\delta^{*}$} & {$\delta=0.5h_{m}$} & {$\delta=h_{m}$} & {$\delta=2h_{m}$} \\ \midrule
    8   & 2E-2  &  7E-3 & 1E-2 &  2E-2 & 7E-3 &  1E-2 \\
    16  & 6E-3 &  4E-4  & 1E-3 &  6E-3    &   5E-4   &  5E-3     \\    
    32  & 3E-3 & 2E-5  &  6E-5 &    4E-3  &  1E-4    &  1E-3     \\
    64  & 2E-4 &  1E-6 &  6E-6 &    7E-4  &    6E-5  &   6E-4    \\ \bottomrule

\end{tabular}
\subcaption{Relative error in the computation of  Stokes single layer potential (with density $\bbf = (x^{2},y^{2},z^{2})$) on  the ellipsoid $\frac{x^{2}}{a^{2}} + \frac{y^{2}}{b^{2}} + \frac{z^{2}}{c^{2}} = 1$, where $a=0.4, b=1, c=1$, for different choices of regularization parameter $\delta$.  Reference solution is computed using $p=64$ spherical harmonics.}
\end{table}

\subsection{Wall clock times for our GPU accelerated code \emph{vs} the spherical harmonics CPU code \label{app:wallclock}}
We report the wall clock time per time step required for our GPU-accelerated code compared to the spherical harmonics CPU code in ~\cref{tab:wall_clock_sph}. We note that the spherical harmonics discretization of degree $p$ contains $2p(p+1)$ points compared to $6(m-1)^2$ for our scheme with $m_{\mathrm{th}}$-order grid. Since, we use four times upsampling, the number of discretization points for evaluating singular quadrature is even higher at $N_{\mathrm{up}} = 6(4m-1)^2$ in our scheme.  Even though our discretization has much higher points than the spherical harmonics, our scheme has lower runtime for higher order grids and hence, scales better than the spherical harmonics as the number of points increase. Since the quadrature scheme requires most of the time in our simulations~\cref{fig:time_dist}, our  scheme can be further accelerated using fast multipole methods(FMMs)~\cite{greengard2008} and can be used to do faster simulations at higher resolutions. We leave the FMM acceleration as the subject for future work.  
\begin{table}
\centering
\caption{Wall clock time per time step (denoted by $T_{\mathrm{step}}$ in \emph{seconds}) for our code with 4x upsampling (on the left) vs the spherical harmonics scheme (on the right). $N$ denotes the number of discretization points. $N_{\mathrm{up}}$ is the upsampled number of discretization points used for the singular quadrature in our scheme. \label{tab:wall_clock_sph}}
\begin{tabular}{SSSS|SSS} \toprule
    {$m$}  & {$N$  } & {$N_{\mathrm{up}}$} & {$T_{\mathrm{step}}$} & {$p$} & {$N$} &  {$T_{\mathrm{step}}$} \\ \midrule
    8   & 294  & 4704 & 0.61 & 8 & 144  &  0.75  \\
    16  & 1350 & 21600 &1.26 & 16  & 544 & 1.92   \\    
    32  & 5766 & 92256 &10.10 & 32  & 2112 & 10.05   \\
    48  & 13254 & 381024 &46.46 & 64 & 8320 & 87.86   \\ \bottomrule

\end{tabular}
\end{table}

\FloatBarrier
\bibliographystyle{abbrv}
\bibliography{master}

\end{document}